\documentclass[notitlepage,leqno,10pt]{article}
\usepackage{latexsym}
\usepackage{amsmath}
\usepackage{amsfonts}
\usepackage{amssymb}

\renewcommand{\d}{\delta }
\newcommand{\D }{\Delta }

\newcommand{\n }{\nabla }

\newcommand{\Sig }{\Sigma}

\newcommand{\ov}{\overline}
\newcommand{\intbar}{\mathop{\int\makebox(-13.5,0){\rule[4pt]{.7em}{0.3pt}}%
\kern-6pt}\nolimits}

\newcommand{\be}{\begin{equation}}
\newcommand{\ee}{\end{equation}}
\newcommand{\bes}{\begin{equation*}}
\newcommand{\ees}{\end{equation*}}
\newcommand{\ba}{\begin{eqnarray}}
\newcommand{\ea}{\end{eqnarray}}
\newcommand{\bas}{\begin{eqnarray*}}
\newcommand{\eas}{\end{eqnarray*}}
\newenvironment{pf}{\noindent{\sc Proof}.\enspace}{\rule{2mm}{2mm}\medskip}

\newcommand{\R}{\mathbb{R}}

\newcommand{\N}{\mathbb{N}}
\def\dint{\displaystyle{\int}}
\author{ Cheikh Birahim NDIAYE}

\date{}

\title{\bf Curvature flows on four manifolds with boundary}
\begin{document}

\newtheorem{lem}{Lemma}[section]
\newtheorem{pro}[lem]{Proposition}
\newtheorem{thm}[lem]{Theorem}
\newtheorem{rem}[lem]{Remark}
\newtheorem{cor}[lem]{Corollary}
\newtheorem{df}[lem]{Definition}

\maketitle

\begin{center}

{\small  SISSA, via Beirut 2-4, 34014 Trieste, Italy.}

\end{center}

\

\

\begin{center} 
{\bf Abstract}
\end{center}

Given a compact four dimensional smooth  Riemannian manifold  \;$(M,g)$\;with smooth boundary, we consider the evolution equation by \;$Q$-curvature in the interior keeping the \;$T$-curvature and the mean curvature to be zero and the evolution equation by \;$T$-curvature at the boundary with the condition that the \;$Q$-curvature  and the mean curvature vanish. Using integral method, we prove global existence and  convergence for the \;$Q$-curvature flow (resp\;$T$-curvature flow) to smooth metric of prescribed \;$Q$-curvature (resp \;$T$-curvature) under conformally invariant assumptions.
\begin{center}

\bigskip

\centerline{\bf AMS subject classification: 35B33, 35J35, 53A30, 53C21}

\end{center}
\section{Introduction}
In the last decades there has been an intensive study on Geometric Flows in order to understand geometrical, analytical, topological and physical problems. Such flows include the mean, inverse mean, Gauss curvature and Willmore flows of submanifolds, Ricci, K\"ahler-Ricci and Calabi flows of manifolds, Yang-Mills flows of connection, Yamabe ,\;$Q$-curvature and other conformal flows of metrics. Just to mention some applications of the the study of such a flow, we have the Geometrization Conjecture, the Riemnann Penrose inequality in General Relativity, Yang-Mills connection, Harmonic maps, some Uniformizations Theorems in K\"ahler geometry and Conformal geometry, ect...\\
\indent
In the context of \;$2$-dimensional geometry, we have the example of the Gauss curvature flow
\begin{equation*}
\frac{\partial g(t)}{\partial t}=-2(K_{g(t)}-\ov{K_{g(t)}})g(t)\;\;\text{on}\;\;\Sig;\\
\end{equation*}
where \;$\Sig$\;is a compact closed surface, \;$K_{g(t)}$\;the Gauss curvature of the evolving Riemannian manifold \;$(M, g(t))$\;and $\ov {K_{g(t)}}$\;the mean value of \;$K_{g(t)}$.
From a Theorem of R. Hamilton\cite{ha} and an other one of B. Chow\cite{ch}, it is proved that for every closed surface \;$\Sig$\; and every initial metric \;$g_0$\;on \;$\Sig$, the initial value problem
\begin{equation*}
\left\{
\begin{split}
&\frac{\partial g(t)}{\partial t}=-2(K_{g(t)}-\ov{K_{g(t)}})g(t)\;\;&\text{on}\;\;\Sig;\\
&g(0)=g_{0}\;\;\;&\text{on}\;\;\Sig.\\
\end{split}
\right.
\end{equation*}
admits a unique globally defined solution which converges exponentially (as $t\rightarrow+\infty$) to a metric (conformal to \;$g_0$) of constant Gauss curvature. Hence providing an other proof of the Classical {\em Uniformization} Theorem for compact closed Riemannian surfaces.

\indent
In \;$4$-dimensional geometry, we have the \;$Q$-curvature flow
\begin{equation*}
\frac{\partial g(t)}{\partial t}=-2(Q_{g(t)}-\ov{Q_{g(t)}})g(t)\;\;\;\;\text{on}\;\;M
\end{equation*}
where  \;$M$\;is a \;$4$-dimensional compact closed manifold,\;$Q_{g(t)}=-\frac{1}{12}(\D_{g(t)} R_{g(t)}-R_{g(t)}^2+3|Ric_{g(t)}|^2)$\;is the \;$Q$-curvature of the evolving Riemannian manifold \;$(M, g(t))$\; and \;$\ov{Q_{g(t)}}$\;its mean value with \;$R_{g(t)}$\;denoting the scalar curvature and \;$Ric_{g(t)}$\;the Ricci tensor.\\
Just to mention an existence and convergence results for this flow we cite the one of S. Brendle see\cite{bren}. He proved that for every four dimensional compact closed manifold \;$M$\;and for every initial metric \;$g_0$\;on \;$M$, if the Paneitz operator \;$P^4_{g_0}$\;(for the definition see \;$\eqref{eq:paneitz}$) is non-negative with trivial kernel and \;$\int_{M}Q_gdV_g<8\pi^2$\;then
\;the initial value problem
\begin{equation*}
\left\{
\begin{split}
&\frac{\partial g(t)}{\partial t}=-2(Q_{g(t)}-\ov{Q_{g(t)}})g(t)\;\;\;\;&\text{on}\;\;M,\\
&g(0)=g_0\;\;\;&\text{on}\;\;M.
\end{split}
\right.
\end{equation*}
 admits a unique globally defined solution which converges (as\;$t\rightarrow+\infty$) to a smooth metric (conformal to \;$g_0$) with constant \;$Q$-curvature.\\
\indent
Still in the context of \;$2$-dimensional geometry, but for compact surfaces with boundary \;$\Sig$, S. Brendle has considered the following two flows:
\begin{equation}\label{eq:gflowb}
\left\{
\begin{split}
&\frac{\partial g(t)}{\partial t}=-2(K_{g(t)}-\ov{K_{g(t)}})g(t)\;\;&\text{on}\;\;\Sig;\\
&k_{g(t)}=0\;\;\;&\text{on}\;\;\partial \Sig.\\
\end{split}
\right.
\end{equation}
\begin{equation}\label{eq:geflowb}
\left\{
\begin{split}
&\frac{\partial g(t)}{\partial t}=-2(k_{g(t)}-\ov{k_{g(t)}})g(t)\;\;&\text{on}\;\;\partial \Sig;\\
&K_{g(t)}=0\;\;\;&\text{on}\;\; \Sig.\\
\end{split}
\right.
\end{equation}
where \;$K_{g(t)}$\;denotes the Gauss curvature of the evolving Riemannian surface \;$(M, g(t))$\; and \;$k_{g(t)}$\;its geodesic curvature. As a result of his study, see \cite{br1}, he proved that if the initial metric has vanishing geodesic curvature at the boundary then the initial value problem corresponding to\;$\eqref{eq:gflowb}$\;has a unique globally define solution which converges to a metric (conformal to the initial one) with constant Gauss curvature and vanishing geodesic curvature. Analogously he show also that if the initial metric has vanishing Gauss curvature in the interior then the initial value problem corresponding to \;$\eqref{eq:geflowb}$\; admits a unique globally defined solution which converges to a metric (conformal to the initial one) with constant geodesic curvature and vanishing Gauss curvature.\\
\vspace{2pt}

\indent
In this paper we are interested to investiguate natural counterparts in \;$4$-dimensional conformal geometry (for manifolds with boundary) of the flows on surfaces with boundary ($\eqref{eq:gflowb}$,\;$\eqref{eq:geflowb}$) considrered by S. Brendle .\\

In the theory of surfaces, it is well known  that the Laplace-Beltrami operator on compact surfaces  with boundary\;$\Sigma$\; and the Neumann operator on the boundary are conformally invariant ones and  govern the transformation laws of the Gauss curvature and the geodesic curvature. In fact under the conformal change of metric \;$g_u=e^{2u}g$,\;we have 
\begin{equation*}
\left\{
     \begin{split}
\D_{g_u}=e^{-2u}\D_g;\\
\frac{\partial }{\partial n_{g_u}}=e^{-u}\frac{\partial }{\partial n_{g}};
     \end{split}
     \right.
\qquad \mbox{and} \qquad
\left\{
     \begin{split}
      -\D_{g}u+K_{g}=K_{g_u}e^{2u}\;\;\text{in}\;\;\Sigma;\\
        \frac{\partial u}{\partial n_{g}}+k_g=k_{g_u}e^{u}\;\;\;\text{on}\;\;\partial \Sigma.
     \end{split}
   \right.
\end{equation*}
where \;$\Delta_{g}$\;(resp. \;$\D_{g_u}$)\;is the Laplace-Beltrami operator of (\;$\Sigma,g$)\;(resp. ($\Sigma,g_u$))\; and\;$K_{g}$\;(resp. \;$K_{g_u}$)\;is the Gauss curvature of ($\Sigma,g$) (resp. of ($\Sigma,g_u$)), \;$\frac{\partial }{\partial n}_{g}$\;(resp \;$\frac{\partial }{\partial n}_{g_u}$)\;is the Neumann operator of (\;$\Sigma ,g$)\;(resp. of ($\Sigma_,g_u$)) and \;$k_{g}$\; (resp.\;$k_{g_u}$) is the geodesic curvature of \;($\partial \Sigma,g$) (resp of ($\partial \Sigma_, g_u$)) .\\
Moreover we have the Gauss-Bonnet formula which relates  \;$\int_{\Sigma}K_{g}dV_{g}+\int_{\partial \Sigma}k_gdS_g$\;and the topology of \;$\Sigma$
\begin{equation}\label{eq:gb}
\int_{\Sigma}K_{g}dV_{g}+\int_{\partial \Sigma}k_gdS_g=2\pi\chi(\Sigma),
\end{equation}
where\;$\chi(\Sigma)$\; is the Euler-Poincar\'e  characteristic of \;$\Sigma$,\;$dV_g$\;is the element area of \;$\Sig$ and \;$dS_g$\;is the line element of \;$\partial \Sig$. Thus\;$ \int_{\Sigma}K_{g}dV_{g}+\int_{\partial \Sigma}k_gdS_g$\;is a topological invariant, hence a conformal one.\\
\indent

In \;$4$-dimensional geometry, we have the Paneitz operator\;$P^4_g$\; and the Q-curvature \;$Q_g$\;defined as follows.
\begin{equation}\label{eq:paneitz}
P^4_g\varphi=\D_{g}^{2}\varphi+div_{g}\left((\frac{2}{3}R_{g}g-2Ric_{g})d\varphi\right);\;\;\;\;\;\;\;Q_g=-\frac{1}{12}(\D_{g}R_{g}-R_{g}^{2}+3|Ric_{g}|^{2})
\end{equation} 
where\;$\varphi$\; is any smooth function on \;$M$ and \;$(M,g)$\;is a \;$4$-dimensional Riemannian manifolds with boundary. We have also that \;$P^4_g$\;is conformally invariant in the following sens
\begin{equation*}
P^{4}_{e^{4w}g}=e^{-4w}P^4_g.
\end{equation*}
Likewise Chang and Qing, see \cite{cq1}, have  discovered a boundary operator \;$P^3_g$\;defined on the boundary of compact four dimensional smooth manifolds and a natural third-order curvature \;$T_g$\;associated to \;$P^3_g$ as follows
\begin{equation*}
P^3_g\varphi=\frac{1}{2}\frac{\partial {\D_g\varphi}}{\partial n_g}+\D_{\hat g}\frac{\partial \varphi}{\partial n_g}-\frac{4}{3}H_g\D_{\hat g}\varphi+(L_g)_{ab}(\nabla_{\hat g})_{a}(\nabla_{\hat g})_{b}+\frac{2}{3}\nabla_{\hat g}H_g.\nabla_{\hat g}\varphi+(F-\frac{R_g}{3})\frac{\partial \varphi}{\partial n_g}.
\end{equation*}
\begin{equation*}
T_g=-\frac{1}{12}\frac{\partial R_g}{\partial n_g}+\frac{1}{2}R_gH_g-<G_g,L_g>+3H_g^3-\frac{1}{3}Tr(L^3)+\D_{\hat g}H_g,
\end{equation*} 
where\;$\varphi$\; is any smooth function on \;$M$,\;\;\;$\hat g$\;is the metric induced by \;$g$\;on\;\;$\partial M$, \;$L_g=(L_g)_{ab}=-\frac{1}{2}\frac{\partial g_{ab}}{\partial n_g}$\\is the second fundamental form of \;$\partial M$,\;\;$H_g=\frac{1}{3}tr(L_g)=\frac{1}{2}g^{ab}L^{ab}$\;(where \;$g^{a,b}$\;are the entries of the inverse \;$g^{-1}$\;of the metric\; $g$)\;is the mean curvature of \;$\partial M$,\;$R^k_{bcd}$\;is the Riemann curvature tensor\;\; $F=R^{a}_{nan}$,\;\;$R_{abcd}=g_{ak}R^{k}_{bcd}$\;(where\;$g_{a,k}$\;are the entries of the metric \;$g$)\;and\;\;$<G_g,L_g>=R_{anbn}(L_g)_{ab}$.\\
\vspace{2pt}

\indent
On the other hand, as the Laplace-Beltrami operator and the Neumann operator govern the transformation laws of the Gauss curvature and the geodesic curvature on compact surfaces with boundary under conformal change of metrics, we have that\;$(P^4_g,P^3_g)$\;does the same for \;$(Q_g,T_g)$\;on compact four dimensional smooth manifolds with boundary. In fact, after a conformal change of metric \;$ g_u=e^{2u}g$\;we have that
\begin{equation}\label{eq:conftlaw}
\left\{
     \begin{split}
P^4_{g_u}=e^{-4u}P^4_g;\\
P^3_{g_u}=e^{-3u}P^3_{g};
     \end{split}
   \right.
\qquad \mbox{and}\qquad
\left\{
\begin{split}
P^4_g+2Q_g=2Q_{ g_u}e^{4u}\;\;\text{in }\;\;M\\
P^3_g+T_g=T_{ g_u}e^{3u}\;\;\text{on}\;\;\partial M.
\end{split}
\right.
\end{equation}
Apart from this analogy we have also an extension of the Gauss-Bonnet formula \;$\eqref{eq:gb}$\; which is known as the Gauss-Bonnet-Chern formula
\begin{equation}\label{eq:gbc}
\int_{M}(Q_{g}+\frac{|W_{g}|^{2}}{8})dV_{g}+\int_{\partial M}(T+Z)dS_g=4\pi^{2}\chi(M)
\end{equation}
where \;$W_g$\;denote the Weyl tensor of \;$(M,g)$\; and \;$ZdS_g$\;is a pointwise conformal invariant. Moreover, it turns out that \;$Z$\;vanishes when the boundary is totally geodesic (by totally geodesic we mean that the boundary \;$\partial M$\;is umbilic and minimal).\\
And setting
\begin{equation*}
\kappa_{P^4_g}=\int_{M}Q_gdV_g;
\end{equation*}
\begin{equation*}
\kappa_{P^3_g}=\int_{\partial M}T_gdS_{g};
\end{equation*}
we have that thanks to \;$\eqref{eq:gbc}$, and to the fact that \;$W_gdV_g$\;and \;$LdS_g$\;are pointwise conformally invariant,\;$\kappa_{P^4_g}+\kappa_{P^3_g}$\;is conformally invariant and will be denoted by
\begin{equation}\label{eq:invar}
\kappa_{(P^4,P^3)}=\kappa_{P^4_g}+\kappa_{P^3_g}.
\end{equation}
Therefore in the context of conformal geometry, the pair \;$(Q_g,T_g)$ is a natural generalization of \;$(K_g,k_g)$. Thus the following flows are also natural generalization to the flows \;$\eqref{eq:gflowb}$,\;$\eqref{eq:geflowb}$\; considered by S. Brendle
\begin{equation*}
\left\{
\begin{split}
&\frac{\partial g(t)}{\partial t}=-2(Q_{g(t)}-\ov{Q_{g(t)}})g(t)\;\;\;&\text {on}\;\;M;\\
&T_{g(t)}=0\;\;&\text{on}\;\;\partial M.
\end{split}
\right.
\end{equation*}
and
\begin{equation*}
\left\{
\begin{split}
&\frac{\partial g(t)}{\partial t}=-2(T_{g(t)}-\ov{T_{g(t)}})g(t)\;\;&\text{on}\;\;\partial M;\\
&Q_{g(t)}=0\;\;&\text{on} \;\;M.
\end{split}
\right.
\end{equation*}
However, since both of the flows are fourth order ones with a third order boundary condition, then from the PDE's point of view it is natural to impose a second boundary condition which is of first order. In this spirit of solvability we will impose that the evolving Riemannian manifold \;$(M, g(t))$\; to be minimal, namely \;$H_{g(t)}=0$. So the natural generalization of S. Brendle flows becomes

\begin{equation}\label{eq:qflowb}
\left\{
\begin{split}
&\frac{\partial g(t)}{\partial t}=-2(Q_{g(t)}-\ov{Q_{g(t)}})g(t)\;\;\;&\text {on}\;\;M;\\
&T_{g(t)}=0\;\;&\text{on}\;\;\partial M;\\
&H_{g(t)}=0&\text{on}\;\;\partial M.
\end{split}
\right.
\end{equation}
and
\begin{equation}\label{eq:tflowb}
\left\{
\begin{split}
&\frac{\partial g(t)}{\partial t}=-2(T_{g(t)}-\ov{T_{g(t)}})g(t)\;\;&\text{on}\;\;\partial M;\\
&Q_{g(t)}=0\;\;&\text{on} \;\;M\\
&H_{g(t)}=0&\text{on}\;\;\partial M.
\end{split}
\right.
\end{equation}
Here in this paper, we consider the following two flows
\begin{equation}\label{eq:qflowbf}
\left\{
\begin{split}
&\frac{\partial g(t)}{\partial t}=-2(Q_{g(t)}-\frac{\ov{Q_{g(t)}}}{\ov F}F)g(t)\;\;\;&\text {on}\;\;M;\\
&T_{g(t)}=0\;\;&\text{on}\;\;\partial M;\\
&H_{g(t)}=0&\text{on}\;\;\partial M.
\end{split}
\right.
\end{equation}
and
\begin{equation}\label{eq:tflowbs}
\left\{
\begin{split}
&\frac{\partial g(t)}{\partial t}=-2(T_{g(t)}-\frac{\ov{T_{g(t)}}}{\ov S}S)g(t)\;\;&\text{on}\;\;\partial M;\\
&Q_{g(t)}=0\;\;&\text{on} \;\;M\\
&H_{g(t)}=0&\text{on}\;\;\partial M.
\end{split}
\right.
\end{equation}
where \;$F$\;is a positive smooth function on \;$M$, and \;$S$\;a positive smooth function on \;$\partial M$. We point out that, since the beginning the {\em bar} means the mean value with respect to \;$g(t)$.\\

\noindent
Definding \;$P^{4,3}_g$\;as follows
\begin{equation}
\left<P^{4,3}_gu,v\right>_{L^2(M)}=\int_{M}\D_g u\D_gvdV_g+\frac{2}{3}R_g\nabla_g u\nabla_g vdV_g-2\int_{M}Ric_g(\nabla_g u,\nabla_g v)dV_g-2\int_{\partial M}L_g( \nabla_{\hat g} u, \nabla_{\hat g} v)dS_g,
\end{equation}
we obtain the following theorems:
\begin{thm}\label{eq:tqf}
Let \;$(M,g_0)$\;be a smooth compact Riemannian manifold with smooth boundary such that \;$T_{g_0}=0$\;and \;$H_{g_0}=0$. Assume\;$P^{4,3}_{g_0}$\;is non-negative,\;$ker P^{4,3}_{g_{0}} \simeq\R$\;and \;$\kappa_{(P^4,P^3)}<4\pi^2$, then the initial boundary value problem corresponding to\;$\eqref{eq:qflowbf}$\;with initial data\;$g_0$\;has a unique globally defined solution which converges to a smooth metric (conformal to\;$g_0$) with \;$Q$-curvature \;$Q_{\infty}$\;verifying\;$Q_{\infty}=\frac{\ov{Q_{\infty}}}{\ov F}F$\; and vanishing \;$T$-curvature and mean curvature.
\end{thm}
\begin{thm}\label{eq:ttf}
Let \;$(M,g_0)$\;be a smooth compact Riemannian manifold with smooth boundary such that \;$Q_{g_0}=0$\;and \;$H_{g_0}=0$. Assume\;$P^{4,3}_{g_0}$\;is non-negative,\;$ker P^{4,3}_{g_{0}}\simeq \R$\;and \;$\kappa_{(P^4,P^3)}<4\pi^2$, then the initial boundary value problem corresponding to\;$\eqref{eq:tflowbs}$\;with initial data\;$g_0$\;has a unique globally defined solution which converges to a smooth metric (conformal to\;$g_0$) with \;$T$-curvature \;$T_{\infty}$\;verifying\;$T_{\infty}=\frac{\ov{T_{\infty}}}{\ov S}S$ and vanishing \;$Q$-curvature and mean curvature.
\end{thm}
\begin{rem}\label{eq:rem}
a) The assumptions in both theorems are conformally invariant.\\
b) Theorem\;$\ref{eq:tqf}$\;and Theorem\;$\ref{eq:ttf}$\;do not cover the case where \;$M=B^4$\;or \;$M=S^4_+$.\\
c) No umbilicity condition for \;$\partial M$\; and no local flatness condition for \;$M$\;is assumed for both theorems.\\
d) We point out that apart the proof of the short-time existence of Theorem\;$\ref{eq:tqf}$, all the other steps remain true with trivial adaptations for Theorem\;$\ref{eq:ttf}$. Hence we will give only a full proof for Theorem\;$\ref{eq:tqf}$, and the proof of the short-time existence for Theorem\;$\ref{eq:ttf}$. 
\end{rem}
Our approach to prove the Theorems above follows closely the one in \cite{bren} and \cite{br1}. However, due to the fact that \;$\partial M\neq\emptyset$, there is an evident difference with the present case, which is mainly in the higher-order estimates of the conformal factor. To obtain the \;$W^{4,2}$-bound, we use the same strategy as in \cite{bren} which is based on deriving differential inequalities for the  \;$W^{4,2}$-norm square of the conformal factor. To face the fact that in our situation the underlying manifold has a boundary we use carefully some Sobolev trace embeddings. An other difference comes in deriving \;$W^{k,2}$-bounds. In fact due to the present of the boundary, to get \;$W^{k,2}$-bound for the conformal factor, it is not {\em easy} to do it by studing differential inequalities for \;$||\D^{\frac{k}{2}}u(t)||^2_{L^2}$\;as it is done in \cite{bren}. However, it is natural to consider \;$||(P^{4,3})^{\frac{k}{2}}u(t)||^2_{L^2}$\;instead. Using this new approach, it is no more possible to use the same arguments as in \cite{bren}, because, there one of the main ingredient was that \;$\D^{\frac{k}{2}}$\;(\;$k$\;{\em even}) beeing a differential operator verifies Leibniz rule, which is not the case for \;$(P^{4,3}_{g_0})^{\frac{k}{2}}$\;(even if \;$k$\;is {\em even}), because beeing a pseudodifferential operator. To overcome the lack of Leibniz rule for $(P^{4,3}_{g_0})^{\frac{k}{2}}$, we use {\em commutator} formula in pseudodifferential calculus.
\vspace{10pt}

\noindent

The structure of the paper is the following. In Section 2 we collect some notations and give some preliminaries like the classical Sobolev embedding theorem and interpolation theorem, a Moser-Trudinger type inequality, and a trace analogue of it envolving the operator \;$P^{4,3}_{g_0}$, and  regularity result. In Section 3 we give the proof of Theorem\;$\ref{eq:tqf}$.  The latter Section is divided into three Subsections. The first one is concerned about the derivation of the evolution equation for the conformal factor, the \;$Q$-curvature and the \;$W^{2,2}$-boundedness of the conformal factor. The second one deals with  the higher-order bound on fixed time interval for the conformal factor. And in the last one we establish the global existence and convergence for the flow. The last Section is devoted to the the proof of the short-time existence for the problem \;$\eqref{eq:tflowbs}$, hence with Remark\;$\ref{eq:rem}$, the proof of Theorem\;$\ref{eq:ttf}$\;follows.
\vspace{10pt}

\noindent
{\bf Acknowledgements:}
The author have been supported by M.U.R.S.T within the PRIN 2006 Variational methods and nonlinear differential equations.

\section{Notations and Preliminaries}
In this brief Section, we collect some useful notations, give two geometric functionals which will play an important role for the derivation of the \;$W^{2,2}$\;bound for the conformal factor. Furthermore, we recall some form of Sobolev embedding and interpolation inequalities since, we will often make use of them . Moreover, we state a Moser-Trudinger type inequality on\;$M$, and derive a trace analogue of it , and state a regularity result.\\

\noindent
In the following, \;$W^{m,p}(M,g)$\;(resp\;$W^{m,p}(\partial M,\hat g)$)\; will stands for the usual Sobolev space of functions on \;$(M,g)$\;(resp \;$(\partial M,\hat g)$) which are of class \;$W^{m,p}$\;in each coordinate system.\\
 We remark that the appeareance of\;$g$\;(resp\;$\hat g$\;) in the definition of the latters sets means that the covariant derivatives are taken with respect to \;$g$\;(resp\;$\hat g$).\\
Large positive constant are always denoted by \;$C$, and the value of \;$C$\;is allowed to vary from formula to formula and also within the same line.\\ 
Given \;$u\in L^1(M,g)$\;(resp\;$L^1(\partial M, \hat g)$,\;$\bar u$\;(resp \;$\bar u_{\partial M}$), denote the mean value of \;$u$\;on \;$M$\;(resp $\partial M$), that is\;$\bar u=(Vol_g(M))^{-1}\int_{M}udV_g$\;(resp \;$\bar u_{\partial M}=(Vol_{\hat g}(\partial M))^{-1}\int_{\partial M}udS_{g}$.\\
Sometimes the subscript in the definition of the mean value of a function defined on\;$\partial M$\;is omitted if there is no possibility of confusion.\\\\
From now on the symbol {\em bar} will means the mean value with respect to the evolving metric unless otherwise stated.\\
The notation \;$\ov{u(t)}_{g_0}$\;means the mean value of \;$u(t)$\; on \;$M$\;with respect to the metric \;$g_0$.\\

As already said, we begin with the definition of two geometric functionals. Setting
\begin{equation}\label{eq;espacef}
H_{\frac{\partial}{\partial n}}=\{u\in H^2(M):\;\;\;\frac{\partial u}{\partial n_{g_0}}=0\};
\end{equation}
\begin{equation}\label{eq:gefq}
\begin{split}
II_{Q,F}(u)=<P^{4,3}_{g_0}u,u>_{L^2(M,g_0)}+4\int_{M}Q_{g_0}udV_{g_0}-
\kappa_{(P^4,P^3)}\log\int_{M}Fe^{4u}dV_g;\;\;\;u\in H_{\frac{\partial }{\partial n}},
\end{split}
\end{equation}
and
\begin{equation}\label{eq:geft}
\begin{split}
II_{T,S}(u)=<P^{4,3}_{g_0}u,u>_{L^2(M,g_0)}+4\int_{M}Q_{g_0}udV_{g_0}-
\frac{4}{3}\kappa_{(P^4,P^3)}\log\int_{M}Se^{4u}dV_g;\;\;\;u\in H_{\frac{\partial }{\partial n}},
\end{split}
\end{equation}

\vspace{4pt}

\indent
We point out that critical points of \;$II_{Q,F}$\;(resp \;$II_{T, S}$) belonging to \;$H_{\frac{\partial}{\partial_n}}$\;give rise to metrics conformal to\;$g_0$\;with \;$Q$-curvature a constant multiple of \;$F$\;(\;$T$-curvature a constant multiple of \;$S$), zero \;$T$-curvature\;(resp zero \;$Q$-curvature) and zero mean curvature. The functional \;$II_{Q,1}$\;(resp \;$II_{T,1}$) has been studied respectively in \cite{nd1} (resp \cite{nd2}).\\

Next, we recall the classical Sobolev embedding theorem. Its proof in the Euclidean case can be found in \cite{fr}, page 28, Theorem 10.2. And the curved version (that we will state) can be obtained by standard covering arguments.
\begin{thm}\label{eq:sobolev}(Sobolev Emdedding)\\
 Let \;$(N,g)$\;be a compact \;$n$-dimensional Riemannian manifold with or without boundary, and let \;$u\in W^{m,r}(N,g)$,\;$1\leq r\leq\infty$. Then for every \;$j\in[0,m[$\;, if \;$\frac{j}{n}+\frac{1}{r}-\frac{m}{n}>0$, then we have
$$
||u||_{W^{j,p}(N,g)}\leq C||u||_{W^{m,r}(N,g)}
$$
where \;$p$\;is given by \;$\frac{1}{p}=\frac{j}{n}+\frac{1}{r}-\frac{m}{n}$, and \;$C$\;depends only on \;$(N,g)$,$m,\; j$\;and\;$r$.
\end{thm}

Now, we recall an interpolation theorem between Sobolev spaces. Its proof in the Euclidean case can be found in \cite{fr}, page 24, Theorem 9.3. And as above, the curved version (that we will state) can be obtained by standard covering arguments.

\begin{thm}\label{eq:interpolation}(Interpolation Inequality)\\
 Let \;$(N,g)$\;be a compact \;$n$-dimensional Riemmanian manifold with or without boundary. Let \;$u\in W^{m,r}(N,g)\cap L^q(N,g)$, \;$1\leq r,q\leq \infty$. Then, we have for any integer \;$j\in[0,m[$, and for any number \;$a\in[\frac{j}{m},1[$, there holds
$$
||D^ju||_{L^p(N,g)}\leq C(||u||_{W^{m,r}(N,g)})^a(||u||_{L^{q}(N,g)})^{1-a},
$$
where \;$p$\;is given by the following expression
$$
\frac{1}{p}=\frac{j}{n}+a(\frac{1}{r}-\frac{m}{n})+(1-a)\frac{1}{q},
$$
and \;$C$\;depends only on $(N,g)$,\;$r,\;q,\;j$\;and \;$a$.
\end{thm}

Next we give two Lemmas whose proof can be found in \cite{nd1}.
\begin{lem}\label{eq:neq}
Suppose \;$ker P^{4,3}_{g_0}\simeq\R$\;and\;$P^{4,3}_{g_0}$\;non-negative then we have that \;$|\cdot|_P$\;is an equivalent norm to\;$||\cdot||_{W^{2,2}}$\;on  \;$\{u\in H_{\frac{\partial }{\partial n}}\;\;\bar u=0\}$
\end{lem}
\begin{pro}\label{eq:mos-tru}
Assume \;$P^{4,3}_{g_0}$\;is a non-negative operator with \;$Ker P^{4,3}_{g_0}\simeq\R$. Then  we have that for all \;$\alpha <16\pi^2$\;there exists a constant \;$C=C(M,g_0,\alpha)$\;such that 
\begin{equation}\label{eq:mts}
\dint_{M}e^{\frac{\alpha(u-\bar u)^2}{\left<P^{4,3}_{g_0}u,u\right>_{L^2(M,g_0)}}}dV_{g_0}\leq C,
\end{equation}
for all\;$u\in H_{\frac{\partial }{\partial n}}$,\;and hence
\begin{equation}\label{eq:mts}
\log\int_{M}e^{4(u-\bar u)}dV_{g_0}\leq C+\frac{4}{\alpha}\left<P^{4,3}_{g_0}u,u\right>_{L^2(M,g_0)}\;\;\forall u\in H_{\frac{\partial }{\partial n}}.
\end{equation}
\end{pro}

Now we give a Proposition which is a trace analogue of the Moser-Trudinger type inequality above.

\begin{pro}\label{eq:mos-trub}
Assume \;$P^{4,3}_{g_0}$\;is a non-negative operator with \;$Ker P^{4,3}_{g_0}\simeq\R$. Then  we have that for all \;$\alpha <12\pi^2$\;there exists a constant \;$C=C(M,g_0,\alpha)$\;such that 
\begin{equation}\label{eq:mts}
\dint_{\partial M}e^{\frac{\alpha(u-\bar u_{\partial M})^2}{\left<P^{4,3}_{g_0}u,u\right>_{L^2(M,g_0)}}}dS_{g_0}\leq C,
\end{equation}
for all\;$u\in H_{\frac{\partial }{\partial n}}$,\;and hence
\begin{equation}\label{eq:mts}
\log\int_{\partial M}e^{3(u-\bar u)}dS_{g_0}\leq C+\frac{9}{4\alpha}\left<P^{4,3}_{g_0}u,u\right>_{L^2(M,g_0)}\;\;\forall u\in H_{\frac{\partial }{\partial n}}.
\end{equation}
\end{pro}
\begin{pf}
First of all, without loss of generality we can assume \;$\bar u_{\partial M}=0$. Following the same argument as in Lemma 2.2 in \cite{cq2}. we get \;$\forall \beta<16\pi^2$\;there exists \;$C=C(\beta,M)$
\begin{equation*}
\dint_{M}e^{\frac{\beta v^2}{\int_{M}|\D_g v|^2dV_{g_0}}}dV_{g_0}\leq C,\;\;\forall v\in H_{\frac{\partial }{\partial n}}\;\;\text{with}\;\;\bar v_{\partial M}=0.
\end{equation*}
From this, using the same reasoning as in Proposition 2.7 in \cite{nd1}, we derive
\begin{equation}\label{eq:mts}
\dint_{M}e^{\frac{\beta v^2}{\left<P^{4,3}_{g_0}v,v\right>_{L^2(M,g_0)}}}dV_{g_0}\leq C,\;\;\forall v\in H_{\frac{\partial }{\partial n}}\;\;\text{with}\;\; \bar v_{\partial M}=0.
\end{equation}
Now let \;$X$\;be a vector field extending the the outward normal at the boundary \;$\partial M$. Using the divergence theorem we obtain
$$
\int_{\partial M}e^{\alpha u^2}dS_{g_0}=\int_{M}div_{g_0}\left(Xe^{\alpha u^2}\right)dV_{g_0}.
$$
Using the formula for the divergence of the product of a vector fied and a function we get
\begin{equation}\label{eq:imp00}
\int_{\partial M}e^{\alpha u^2}dS_{g_0}=\int_{M}\left(div_{g_0}X+2u\alpha \n_{g}u\n_{g_0}X \right)e^{\alpha u^2}dV_{g_0}.
\end{equation}
Now we suppose \;$<P^{4,3}_{g_0}u,u>_{L^2(M,g_0)}\leq 1$, then since the vector field \;$X$\;is smooth we have
\begin{equation}\label{eq:imp0}
\left|\int_{M}div_{g_0}Xe^{\alpha u^2}dV_{g_0}\right|\leq C;
\end{equation}
thansk to \;$\eqref{eq:mts}$.
Next let us show that 
\begin{equation*}
\left|\int_{M}2\alpha u\n_{g_0}u\n_{g_0}X e^{\alpha u^2}dV_{g_0}.
\right|\leq C
\end{equation*}
Let \;$\epsilon>0$\;small and let us set \;$$p_1=\frac{4}{3-\epsilon},\;\;p_2=4,\;\;p_3=\frac{4}{\epsilon}.$$
It is easy to check that 
$$\frac{1}{p_1}+\frac{1}{p_2}+\frac{1}{p_3}=1.$$
Using Young's inequality we obtain
$$
\left|\int_{M}2\alpha u\n_{g_0}u\n_{g_0}X e^{\alpha u^2}dV_{g_0}
\right|\leq C||u||_{L^{\frac{4}{\epsilon}}(M,g_0)}||\n_{g_0}u||_{L^4(M,g_0)}\left(\int_{M}e^{\alpha\frac{4}{3-\epsilon}u^2}dV_{g_0}\right)^{\frac{3-\epsilon}{4}}.
$$
On the other hand, Lemma\;$\eqref{eq:neq}$\; and Sobolev embedding theorem imply 
$$
||u||_{L^{\frac{4}{\epsilon}}(M,g_0)}\leq C;
$$
and
$$
||\n_{g}u||_{L^4(M,g_0)}\leq C.
$$
Furthermore from the fact that \;$\alpha<12\pi^2$, by taking \;$\epsilon$\;sufficiently small and using \;$\eqref{eq:mts}$, we obtain
$$
\left(\int_{M}e^{\alpha\frac{4}{3-\epsilon}u^2}dV_{g_0}\right)^{\frac{3-\epsilon}{4}}.
$$
Thus we arrive to
\begin{equation}\label{eq:imp1}
\left|\int_{M}2\alpha u\n_{g_0}u\n_{g_0}X e^{\alpha u^2}dV_{g_0}
\right|\leq C.
\end{equation}
Hence \;$\eqref{eq:imp00}$,\;$\eqref{eq:imp0}$\;and \;$\eqref{eq:imp1}$\;imply
$$
\int_{\partial M}e^{\alpha u^2}dS_{g_0}\leq C,
$$
as desired. So the first point of the Lemma is proved.\\
Now using the algebraic inequality 
$$3ab\leq 3\gamma^2a^2+\frac{3b^2}{4\gamma^2},$$
we have that the second point follows directly from the first one. Hence the Lemma is proved.
\end{pf}

Now, we give a regularity result whose proof is a trivial adaptation of the arguments in Proposition 2.4 in \cite{nd1}.
\begin{lem}\label{eq:reg}
 Let \;$f\in W^{m,p}(M,g_0)$, and \;$h\in W^{m,p}(\partial M,\hat g_0)$. Assume \;$u\in H_{\frac{\partial }{\partial_n}}$\;is a weak solution to 
\begin{equation*}
\left\{
\begin{split}
P^4_gu&=f\;\;&\text{in}\;\;M;\\
P^3_gu&=h\;\;&\text{on}\;\;\partial M.
\end{split}
\right.
\end{equation*}
then 
$$
||u||_{W^{m,p}(M,g_0)}\leq C(||f||_{W^{m,p}(M,g_0)}+||h||_{W^{m,p}(\partial M,\hat g_0)}).
$$
\end{lem}

\section{Proof of Theorem\;$\ref{eq:tqf}$}
This Section is concerned with the proof of Theorem\;$\ref{eq:tqf}$. It is divided in three Subsections. The first one deals with the short-time existence, the evolution of the \;$Q$-curvature and the derivation of the \;$W^{2,2}$-bound for the conformal factor. The second one is concerned about the higher-order estimates on fixed time interval  for the conformal factor. The last one is about higher-order uniform estimates for the conformal factor and the convergence of the flow.

\subsection{Evolution equation of the conformal factor,\;$Q$-curvature and \;$W^{2,2}$-estimates}
In this Subsection, we show the equivalence between our initial boundary value problem and a scalar quasilinear parabolic initial boundary value problem. From this, we derive the short-time existence of solution. Moreover, we derive the evolution equation for the \;$Q$-curvature and prove also the monotonicity of the geometric functional \;$II_{Q,F}$\;along the flow. Furthermore using the Moser-Trudinger type inequality (see Lemma\;$\ref{eq:mos-tru}$), and the monotonocity of \;$II_{Q,F}(u(t))$, we get the \;$W^{2,2}$-boundedness for the conformal factor\;$u(t)$. \\

Now as said in the introduction of the Subsection, we start by giving an equivalent formulation of our evolution problem which turns out to be a a scalar quasilinear parabolic initial boundary value one, and the derivation of the evolution equation for the \;$Q$-curvature.
\begin{lem}\label{eq:scfor}
The evolution equation \;$\eqref{eq:qflowbf}$\;conserve the conformal structure of \;$M$. Moreover, setting the evolving metric of the initial boundary value problem corresponding to \;$g(t)$\; to be \;$g(t)=e^{2u(t)}g_0$, we have that the conformal factor\;$u(t)$\; satisfy the following quasilinear parabolic boundary value problem.
\begin{equation}\label{eq:evu}
\left\{
\begin{split}
&\frac{\partial u(t)}{\partial t}=-\frac{1}{2}\left(e^{-4u(t)}(P^4_{g_0}u(t)+2Q_{g_0})-2\frac{\ov{Q_{g(t)}}}{\ov F}F\right);\\
&P^3_{g_0}u(t)=0\;\;\text{on}\;\;\partial M;\\
&\frac{\partial u(t)}{\partial n_{g_0}}=0\;\text{on}\;\;\partial M;\\
&u(0)=0\;\;\text{on}\;\;M.
\end{split}
\right.
\end{equation}
Furthermore we have that the volume is preserved by the flow \;$\eqref{eq:qflowbf}$ (hence  \;$\ov{Q_g(t)}$, too ) and \;$(Q_{g(t)}-\frac{\ov{Q_{g(t)}}}{\ov F}F)$\; evolves as follows
\begin{equation}
\frac{\partial  (Q_{g(t)}-\frac{\ov {Q_{g(t)}}}{\ov F}F)}{dt}=4(Q_{g(t)}-\frac{\ov {Q_{g(t)}}}{\ov F})-\frac{1}{2}P^4_{g(t)}(Q_{g(t)}-\frac{\ov {Q_{g(t)}}}{\ov F})-4\frac{\ov{Q_{g(t)}}}{\ov F}F\ov{\frac{F}{\ov F}(Q_{g(t)}-\frac{\ov {Q_{g(t)}}}{\ov F})}.
\end{equation}
\end{lem}
\begin{pf}
The fact that the evolution equation\;$\eqref{eq:qflowbf}$ conserves the conformal structure of \;$M$ is evident.\\
Now writting the evolving metric \;$g(t)=e^{2u(t)}g_0$\; we have that by differentiationg \;$g(t)$\;
\begin{equation}
\frac{\partial g(t)}{\partial t}=2\frac{\partial u(t)}{\partial t}e^{2u(t)}g_0=2\frac{\partial u(t)}{\partial t}g(t).
\end{equation}
On the other hand using the evolution equation \;$\eqref{eq:qflowbf}$\;we get
\begin{equation}\label{eq:ftw1}
\frac{\partial u(t)}{\partial t}=-(Q_{g(t)}-\frac{\ov{Q_{g(t)}}}{\ov F}F)\;\;\text{on}\;\;M.
\end{equation}
Next using the fact that the Paneitz operator governs the transformation law of \;$Q$-curvature (for conformal metrics), see \;$\eqref{eq:conftlaw}$, we have
\begin{equation}
2Q_{g(t)}=e^{-4u(t)}\left(P^4_{g_0}u(t)+2Q_{g_0}\right)\;\;\text{on}\;\;M.
\end{equation}
Hence \;$\eqref{eq:ftw1}$\;becomes
\begin{equation}\label{eq:eequw}
\frac{\partial u(t)}{\partial t}=-\frac{1}{2}\left(e^{-4u(t)}(P^4_{g_0}u(t)+2Q_{g_0})-2\frac{\ov{Q_{g(t)}}}{\ov F}F\right)\;\;\text{on}\;\;M.
\end{equation}
Furthermore using also the fact that the Chang-Qing operator governs the transformation law of \;$T$-curvature (still for conformal metrics), see \;$\eqref{eq:conftlaw}$, we obtain
\begin{equation}
P^3_{g_0}u(t)+T_{g_0}=T_{g(t)}e^{3u(t)}\;\;\text{on}\;\;\partial M.
\end{equation}
Now from the fact that the initial metric \;$g_0$\;and the evolving metric \;$g(t)$\; all are \;$T$-flat, we infer that 
\begin{equation}\label{eq:p3c}
P^3_{g_0}u(t)=0.
\end{equation}
Next using the fact that the Neumann operator governs the transformation law of mean curvature (within a conformal class), we get 
\begin{equation}\label{eq:nc}
\frac{\partial u}{\partial n}_{g_0}+H_{g_0}=H_{g}e^{u}\;\;\text{on}\;\;\partial M.
\end{equation}
Hence using the fact that the initial metric and the evolving one all are such that the underlying manifold endowed with is minimal, we derive that
\begin{equation}
\frac{\partial u(t)}{\partial n_{g_0}}=0.
\end{equation}
So \;$\eqref{eq:eequw}$,\;$\eqref{eq:p3c}$,\;$\eqref{eq:nc}$\; and the fact that we start we the initial metric \;$g_0$\;imply that the first point of the Lemma is proved.\\
Now let show that the volume is invariant with respect to time ($t$). First of all it is well known that 
\begin{equation}\label{eq:vol}
dV_{g(t)}=e^{4u(t)}dV_{g_0}.
\end{equation}
Hence from the definition of the volume and \;$\eqref{eq:vol}$\;we obtain 
\begin{equation}\label{eq:forv}
Vol_{g(t)}(M)=\int_{M}e^{4u(t)}dV_{g_0}.
\end{equation}
Now differentiating with respect to time \;$\eqref{eq:forv}$\;and using again \;$\eqref{eq:vol}$, we get
\begin{equation}
\frac{d Vol_{g(t)}(M)}{dt}=4\int_{M}\frac{\partial u(t)}{\partial t}dV_{g(t)}.
\end{equation}
Thus using \;$\eqref{eq:ftw1}$\;we arrive to
\begin{equation*}
\frac{d Vol_{g(t)}(M)}{dt}=4\int_{M}(F\frac{\ov{ Q_{g(t)}}}{\ov F}-Q_{g(t)})dV_{g(t)}=0.
\end{equation*}
Hence the volume is invariant (with respect to time). \\
On the other hand since \;$T_{g(t)}=0$\; then from the general fact that the total integral \;$(Q,T)$-curvature is conformally invariant, we get\;$\int_{M}Q_{g(t)}dV_{g(t)}$\;is invariant with time. Thus from the invariance of the volume, we infer that \;$\ov{Q_{g(t)}}$\;is also invariant with respect to time.\\
Now let us derive the evolution equation of \;$(Q_{g(t)}-\frac{\ov{Q_{g(t)}}}{\ov F}F)$. 
Firts recalling that \;$\ov{Q_{g(t)}}$\;is invariant with respect to time we have 
\begin{equation}\label{eq:der0}
\frac{\partial (Q_{g(t)}-\frac{\ov{Q_{g(t)}}}{\ov F}F)}{\partial t}=\frac{\partial Q_{g(t)}}{\partial t}+\frac{\ov {Q_{g(t)}}}{\ov{F}^2}F\frac{d\ov{F}}{dt}.
\end{equation}
Now using the fact that \;$Q_{g(t)}=\frac{1}{2}e^{-4u(t)}(P^4_{g_0}u(t)+2Q_{g_0})$, see \;$\eqref{eq:conftlaw}$\; and Leibniz rule we get
\begin{equation}\label{eq:der1}
\frac{\partial Q_{g(t)}}{\partial t}=-2\frac{\partial u(t)}{\partial t}e^{-4u(t)}(P^4_{g_0}u(t)+2Q_{g_0})+\frac{1}{2}P^4_{g(t)}(\frac{\partial u(t)}{\partial t}).
\end{equation}
Thus using again the formula for \;$Q_{g(t)}$\; we have that \;$\eqref{eq:der1}$\;becomes
\begin{equation}\label{eq:der2}
\frac{\partial Q_{g(t)}}{\partial t}=-4\frac{\partial u(t)}{\partial t}Q_{g(t)}+\frac{1}{2}P^4_{g(t)}(\frac{\partial u(t)}{\partial t})
\end{equation}
Now \;$\eqref{eq:ftw1}$ \;and \;$\eqref{eq:der2}$\;give 
\begin{equation}\label{eq:der3}
\frac{\partial Q_{g(t)}}{\partial t}=4(Q_{g(t)}-\frac{\ov{Q_{g(t)}}}{\ov F}F)Q_{g(t)}-\frac{1}{2}P^4_{g(t)}(Q_{g(t)}-\frac{\ov{Q_{g(t)}}}{\ov F}F).
\end{equation}
On the other hand, by trivial calculations, one can easily check that the following holds
\begin{equation}\label{eq:derbf}
\frac{d\ov{F}}{dt}=-4\ov{F(Q_{g(t)}-\frac{\ov{Q_{g(t)}}}{\ov F})}.
 \end{equation}
Thus \;$\eqref{eq:der0}$,\;$\eqref{eq:der3}$, and \;$\eqref{eq:derbf}$\;imply 
\begin{equation}
\frac{\partial  (Q_{g(t)}-\frac{\ov {Q_{g(t)}}}{\ov F}F)}{dt}=4(Q_{g(t)}-\frac{\ov {Q_{g(t)}}}{\ov F})-\frac{1}{2}P^4_{g(t)}(Q_{g(t)}-\frac{\ov {Q_{g(t)}}}{\ov F})-4\frac{\ov{Q_{g(t)}}}{\ov F}F\ov{\frac{F}{\ov F}(Q_{g(t)}-\frac{\ov {Q_{g(t)}}}{\ov F})}.
\end{equation}
So the last point of the Lemma holds. Hence the proof of the Lemma is concluded.
\end{pf}

Next, we give a Proposition which shows the short-time existence of the flow, the monotonicity of \;$II_{Q,F}$\;along it, and the \;$W^{2,2}$\;bound of the conformal factor.
\begin{pro}\label{eq:eh2b}
Under the assumptions of Theorem\;$\ref{eq:tqf}$, we have that the initial boundary value problem corresponding to \;$\eqref{eq:qflowbf}$\;has a unique short time solution \;$g(t)=e^{2u(t)}g_0$. Moreover \;the conformal factor\;$u(t)$\;satisfies 
\begin{equation}\label{eq:dec}
\frac{d II_{Q,F}(u(t))}{dt}=-4\int_{M}\left(Q_{g(t)}-\frac{\ov {Q_{g(t)}}}{\ov F}F\right)^2dV_{g(t)}.
\end{equation}
Furthermore there holds
\begin{equation}\label{eq:h2b}
||u(t)||_{W^{2,2}(M,g_0)}\leq C \;\;\;\;\;\forall t.
\end{equation}
\end{pro}
\begin{pf}
By Lemma\;$\ref{eq:scfor}$, we have that the evolution equation in consideration is equivalent to the following quasilinear parabolic BVP
\begin{equation}\label{eq:eprob}
\left\{
\begin{split}
&\frac{\partial u(t)}{\partial t}=-\frac{1}{2}\left(e^{-4u(t)}(P^4_{g_0}u+2Q_{g_0})-2\frac{\ov{Q_{g(t)}}}{\ov F}F\right);\\
&P^3_{g_0}u(t)=0\;\;\text{on}\;\;\partial M;\\
&\frac{\partial u(t)}{\partial g_0}=0\;\text{on}\;\;\partial M;\\
&u(0)=0\;\;\text{on}\;\;M.
\end{split}
\right.
\end{equation}
Since by asummption \;$P^{4,3}_{g_0}$\;is non-negative with trivial kernel, then the problem \;$\eqref{eq:eprob}$\; is parabolic. Hence the theory for short time existence for scalar parabolic evolution equation ensure the existence of a unique short-time solution.\\
Now let us show that 
\begin{equation}\label{eq:der1f}
\frac{d II_{Q,F}(u(t))}{dt}=4\int_{M}(Q_{g(t)}-\frac{\ov{Q_{g(t)}}}{\ov F}F)\frac{\partial u(t)}{\partial t}dV_{g(t)}.
\end{equation}
From the fact that \;$P^{4,3}_{g_0}$\;is self-adjoint, we infer that
\begin{equation}\label{eq:der2f}
\frac{d II_{Q,F}(u(t))}{dt}=2<P^{4,3}_{g_0}u(t),\frac{\partial u(t)}{\partial t}>_{L^2(M,g_0)}+4\int_{M}Q_{g_0}\frac{\partial u(t)}{\partial t}dV_{g_0}-4\kappa_{(P^4,P^3)}\frac{\int_{M}Fe^{4u(t)}\frac{\partial u(t)}{\partial t}dV_{g_0}}{\int_{M}Fe^{4u(t)}dV_{g_0}}.
\end{equation}
Now recalling that \;$dV_{g(t)}=e^{4u(t)}dV_{g_0}$\;we have that \;$\eqref{eq:der2f}$\;becomes
\begin{equation}
\frac{d II_{Q,F}(u(t))}{dt}=2<P^{4,3}_{g_0}u(t),\frac{\partial u(t)}{\partial t}>_{L^2(M,g_0)}+4\int_{M}Q_{g_0}\frac{\partial u(t)}{\partial t}dV_{g_0}-\int_{M}4\frac{\ov{Q_{g(t)}}}{\ov F}F\frac{\partial u(t)}{\partial t}dV_{g(t)}.
\end{equation}
Next, from \;$\eqref{eq:conftlaw}$\;and Lemma\;$\ref{eq:scfor}$, we get
\begin{equation}\label{eq:bvput}
\left\{
\begin{split}
P^4_{g_0}u(t)+2Q_{g_0}&=2 Q_{g(t)}e^{4u(t)}\;\;&\text{in}\;\;M;\\
P^3_gu(t)&=0\;\;&\text{on}\;\;\partial M;\\
\frac{\partial u(t)}{\partial n_{g_0}}&=0\;\;&\text{on}\;\;\partial M.
\end{split}
\right.
\end{equation}
Using \;$\eqref{eq:bvput}$\;we derive
\begin{equation}\label{eq:intbvp}
2<P^{4,3}_{g_0}u(t),\frac{\partial u(t)}{\partial t}>_{L^2(M,g_0)}+4\int_{M}Q_{g_0}\frac{\partial u(t)}{\partial t}dV_{g_0}=4\int_{M}Q_{g(t)}e^{4u(t)}\frac{\partial u(t)}{\partial t}dV_{g_0}.
\end{equation}
Furthermore using \;$\eqref{eq:vol}$, we have that \;$\eqref{eq:intbvp}$\;becomes
\begin{equation}
2<P^{4,3}_{g_0}u,\frac{\partial u}{\partial t}>_{L^2(M,g_0)}+4\int_{M}Q_{g_0}\frac{\partial u}{\partial t}dV_{g_0}=4\int_{M}Q_{g(t)}\frac{\partial u}{\partial t}dV_{g(t)}.
\end{equation}
Thus we obtain
\begin{equation}
\frac{d II(u(t))}{dt}=4\int_{M}(Q_{g(t)}-\frac{\ov{Q_{g(t)}}}{\ov F}F)\frac{\partial u(t)}{\partial t}dV_{g(t)}.
\end{equation}
Hence the claim \;$\eqref{eq:der1f}$\;is proved.
Now recalling that\;$\frac{\partial u(t)}{\partial t}=-(Q_{g(t)}-\frac{\ov{Q_{g(t)}}}{\ov F}F)$\;(see\;$\eqref{eq:ftw1}$ ) we have that \;$\eqref{eq:der1f}$\;becomes
\begin{equation*}
\frac{d II_{Q,F}(u(t))}{dt}=-4\int_{M}(Q_{g(t)}-\frac{\ov{Q_{g(t)}}}{\ov F}F)^2dV_{g(t)}.
\end{equation*}
as desired.\\
Next let us show that
\begin{equation}
||u(t)||_{W^{2,2}(M,g_0)}\leq C\;\;\;\forall t.
\end{equation}
First of all using\;$\eqref{eq:dec}$\;we have that the Energy functional is decreasing along the the flow, hence we infer that
\begin{equation}\label{eq:upp}
II_{Q,F}(u(t))\leq C.
\end{equation}
Now suppose \;$\kappa_{(P^4,P^3)}<0$. We have by Jensen's inequality that 
\begin{equation}
II_{Q,F}(u(t))\geq <P^{4,3}_{g_0}u(t),u(t)>_{L^2(M,g_0)}+4\int_{M}Q_{g_0}(u(t)-\ov {u(t)}_{g_0})dV_{g_0}-C.
\end{equation}
Next using H\"{o}lder inequality, Poincar\'{e} inequality, and Lemma\;$\ref{eq:neq}$\; we get
\begin{equation}\label{eq:lowb1}
II_{Q,F}(u(t))\geq \beta||u(t)-\ov{ u(t)}_{g_0}||_{W^{2,2}(M,g_0)}-C.
\end{equation}
for some \;$\beta>0$.\\
Now if \;$0<\kappa_{(P^4,P^3)}<4\pi^2$, we have by using Lemma\;$\ref{eq:neq}$\; Poincare inequality and Moser-Trudinger type inequality (see Lemma \;$\ref{eq:mos-tru}$), we obtain
\begin{equation}\label{eq:lowb2}
II_{Q,F}(u(t))\geq \d||u(t)-\ov {u(t)}_{g_0}||_{W^{2,2}(M,g_0)}-C
\end{equation}
for some \;$\delta>0$. Thus \;$\eqref{eq:lowb1}$\;and \;$\eqref{eq:lowb2}$\;imply that in both cases, there holds
\begin{equation}\label{eq:lowb3}
II_{Q,F}(u(t))\geq \gamma||u(t)-\ov {u(t)}_{g_0}||_{W^{2,2}(M,g_0)}-C
\end{equation}
for some \;$\gamma>0$. Hence to prove the \;$W^{2,2}$-boundedness of \;$u(t)$, it is sufficient to prove that
\begin{equation}\label{eq:bbu}
-C\leq\ov {u(t)}_{g_0}\leq C.
\end{equation}
To do this we first use Moser-Trudinger type inequality ( see Lemma \;$\ref{eq:mos-tru}$), \;$\eqref{eq:lowb3}$ and\;$\eqref{eq:upp}$\;  to infer that
\begin{equation}\label{eq:upexp}
\int_{M}e^{4(u(t)-\ov {u(t)}_{g_0})}dV_{g_0}\leq C.
\end{equation}
Next from the fact that the volume is conserved, we have that \;$\eqref{eq:upexp}$\;implies
\begin{equation}\label{eq:lobu}
\ov{ u(t)}_{g_0}\geq -C.
\end{equation}
Furthermore from the conservation of the volume and Jensen's inequality, we derive
\begin{equation}\label{eq:upbu}
\ov{u(t)}_{g_0}\leq C.
\end{equation}
Hence \;$\eqref{eq:lobu}$\;and \;$\eqref{eq:upbu}$\;gives
\begin{equation*}
-C\leq \ov{u(t)}_{g_0}\leq C.
\end{equation*}
Thus  \;$\eqref{eq:bbu}$\;is proved.\\
Now using \;$\eqref{eq:lowb3}$\;and \;$\eqref{eq:bbu}$\;we get 
\begin{equation}
||u(t)||_{W^{2,2}(M,g_0)}\leq C. 
\end{equation}
So the \;$W^{2,2}$-boundedness of \;$u(t)$\;is proved. Hence ending the proof of the Proposition.
\end{pf}
\subsection{Higher order a priopri estimates on fixed time interval}
This Subsection is concerned about higher-order estimates of the conformal factor \;$u(t)$\;on fixed time interval. We start with a Lemma showing the \;$W^{4,2}$-estimate.
\begin{lem}
Suppose that the assumptions of Theorem\;$\ref{eq:tqf}$\;holds and 
let \;$u(t)$\;be the unique short-time solution to initial scalar parabolic boundary value problem  \;$\eqref{eq:evu}$. Then for every \;$T>0$\;such that \;
$u(t)$\;is defined on $[0, T[$, there exists \;$C=C(T)>0$\;such that
\begin{equation}
||u(t)||_{W^{4,2}(M,g_0)}\leq C \;\;\;\;\;\forall t\in [0, T[.
\end{equation}
\end{lem}
\begin{pf}
Let us set
\begin{equation}
v(t)=e^{2u(t)}\frac{\partial u(t)}{\partial t}.
\end{equation}
Then using the evolution equation for \;$u(t)$\;(see \;$\eqref{eq:evu}$) we get
\begin{equation}
v(t)=\frac{1}{2}\left(-e^{-2u(t)}(P^{4}_{g_0}u(t)+2Q_{g_0})-2e^{2u(t)}\frac{\ov{Q_{g(t)}}}{\ov F}F\right).
\end{equation}
From this last formula we derive the expression of \;$P^{4}_{g_0}u(t)$\;in term of \;$v(t)$\; as follows
\begin{equation}\label{eq:exppu}
P^{4}_{g_0}u(t)=-2ve^{2u(t)}-2Q_{g_0}+2e^{4u(t)}\frac{\ov{Q_{g(t)}}}{\ov F}F.
\end{equation}
Now let us show that \;$\left(\int_{M}(P^{4,3}_{g_0}u(t))^2dV_{g_0}\right)$\; satisfies a differential inequality that will allow us get an upper bound for it. Hence by the properties of \;$P^{4,3}_{g_0}$, we obtain the \;$W^{4,2}$-bounds for \;$u(t)$\;as desired.\\ 
Using the rule of differentiation under the sign integral, we get
\begin{equation}\label{eq:der2p}
\frac{d}{d t}\left(\int_{M}(P^{4,3}_{g_0}u(t))^2dV_{g_0}\right)=2\int_{M}P^{4,3}_{g_0}u(t) P^{4,3}_{g_0}(\frac{\partial u(t)}{\partial t})dV_{g_0}.
\end{equation}
Next, using\;$\eqref{eq:exppu}$\;and the fact that \;$P^3_{g_0}u(t)=0$, we have that\;$\eqref{eq:der2p}$\;becomes
\begin{equation}\label{eq:der2p1}
\frac{d}{d t}\left(\int_{M}(P^{4,3}_{g_0}u(t))^2dV_{g_0}\right)=2\int_{M}\left(-2ve^{2u(t)}-2Q_{g_0}+2e^{4u(t)}\frac{\ov{Q_{g(t)}}}{\ov F}F\right)P^{4,3}_{g_0}(\frac{\partial u(t)}{\partial t})dV_{g_0}.
\end{equation}
Using the definition of \;$v(t)$\;and expanding the left hand side of \;$\eqref{eq:der2p1}$, we get
\begin{equation}\label{eq:der3p}
\frac{d}{d t}\left(\int_{M}(P^{4,3}_{g_0}u(t))^2dV_{g_0}\right)=\int_{M}\left(-4ve^{2u(t)}P^{4,3}_{g_0}(ve^{-2u(t)})+(4e^{4u(t)}\frac{\ov{Q_{g(t)}}}{\ov F}F-4Q_{g_0})P^{4,3}_{g_0}(v(t)e^{-2u(t)})\right)dV_{g_0}
\end{equation}
On the other hand, from the definition of the operator \;$P^{4,3}_{g_0}$, we have that \;$\eqref{eq:der3p}$\;becomes
\begin{equation}
\begin{split}
\frac{d}{d t}\left(\int_{M}(P^{4,3}_{g_0}u(t))^2dV_{g_0}\right)=\int_{M}\left[-4\D_{g_0}(v(t)e^{2u(t)})\D_{g_0}(v(t)e^{-2u(t)})-\frac{8}{3}R_{g_0}\n_{g_0}(v(t)e^{2u(t)})\n_{g_0}(v(t)e^{-2u(t)})\right]dV_{g_0}\\+\int_{M}8Ric_{g_0}(\n_{g_0}(v(t)e^{2u(t)}),\n_{g_0}(v(t)e^{-2u(t)}))dV_{g_0}+8\int_{\partial M}L_{g_0}(\n_{\hat g_0}(v(t)e^{2u(t)}),\n_{\hat g_0}(v(t)e^{-2u(t)}))dS_{g_0}\\+\int_{M}\left[(4e^{4u(t)}\frac{\ov{Q_{g(t)}}}{\ov F}F-4Q_{g_0})P^{4,3}_{g_0}(v(t)e^{-2u(t)})\right]dV_{g_0}.
\end{split}
\end{equation}
Now using the identity 
$$
\D_{g_0}(v(t)e^{2u(t)})\D_{g_0}(e^{-2u(t)}v(t))=\left(\D_{g_0}v(t)+4|\n_{g_0}u(t)|^2v(t)\right)^2-\left(4(\n_{g_0}u(t),\n_{g_0}v(t))+2\D_{g_0}u(t)v(t)\right)^2$$
we obtain
\begin{equation}\label{eq:lapar}
\D_{g_0}(v(t)e^{2u(t)})\D_{g_0}(e^{-2u(t)}v(t))-\frac{1}{2}(\D_{g_0}v(t))^2\geq  -C\left(|\n_{g_0}u(t)|^4v(t)^2+(\n_{g_0}u(t)\n_{g_0}v(t))^2+(\D_{g_0}u(t))^2v(t)^2\right).
\end{equation}
On the other hand using  H\"{o}lder inequality we get
\begin{equation}\label{eq:guv}
 \int_{M}|\n_{g_0}u(t)|^4v(t)^2dV_{g_0}\leq ||\n_{g_0}u||_{L^8(M,g_0)}^4||v||_{L^4(M,g_0)}^2 ;
\end{equation}
Now using Sobolev embedding Theorem (see Theorem\;$\ref{eq:sobolev}$), one can easily see that the following hold
$$
||u(t)||_{W^{1,8}(M,g_0)}\leq C||u(t)||_{W^{2,\frac{8}{3}}(M,g_0)}.
$$
On the other hand, since \;$-C\leq \ov{u(t)}_{g_0}\leq C$, then we have 
$$
||u(t)||_{W^{2,\frac{8}{3}}(M,g_0)}\leq C||\n_{g_0}^2u(t)||_{L^{\frac{8}{3}}(M,g_0)}.
$$
Next using interpolation inequality (see Theorem\;$\ref{eq:interpolation}$), we have
$$
||\n_{g_0}^2u(t)||_{L^{\frac{8}{3}}(M,g_0)}\leq C||\n_{g_0}^2u(t)||_{L^2(M,g_0)}^{\frac{3}{4}}||\n_{g_0}^2u(t)||_{W^{2,2}(M,g_0)}^{\frac{1}{4}}.
$$
Thus we arrive to
$$
||\n_{g_0}^2u(t)||_{L^{\frac{8}{3}}(M,g_0)}\leq C||u(t)||_{W^{2,2}(M,g_0)}^{\frac{3}{4}}||u(t)||_{W^{4,2}(M,g_0)}^{\frac{1}{4}}.
$$
Hence we obtain
\begin{equation}\label{eq:18u}
 ||u(t)||_{W^{1,8}}\leq C||u(t)||_{W^{2,2}(M,g_0)}^{\frac{3}{4}}||u(t)||_{W^{4,2}(M,g_0)}^{\frac{1}{4}};
\end{equation}
On the other hand, using again the same interpolation Theorem as above, we get
\begin{equation}\label{eq:04v}
 ||v(t)||_{L^4(M,g_0)}\leq C||v(t)||_{L^2(M,g_0)}^{\frac{1}{2}}||v(t)||_{W^{2,2}(M,g_0)}^{\frac{1}{2}}.
\end{equation}
Hence \;$\eqref{eq:guv}$,\;\;$\eqref{eq:18u}$\;and \;$\eqref{eq:04v}$\;imply
\begin{equation}\label{eq:guv1}
 \int_{M}|\n_{g_0}u(t)|^4v(t)^2dV_{g_0}\leq  C||u(t)||_{W^{2,2}(M,g_0)}^3||u(t)||_{W^{4,2}(M,g_0)}||v(t)||_{L^2(M,g_0)}||v(t)||_{W^{2,2}(M,g_0)};
\end{equation}
Using again H\"{o}lder inequality we get
\begin{equation}\label{eq:nunv}
\int_{M}(\n_{0}u(t),\n_{g_0}v(t))^2dV_{g_0}\leq||\n_{g_0}u(t)||_{L^8(M,g_0)}^2||\n_{g_0}v(t)||_{L^{\frac{8}{3}}(M,g_0)}^2;
\end{equation}
Now using interpolation as above, we get
$$
||\n_{g_0}v(t)||_{L^{\frac{8}{3}}(M,g_0)}\leq C||\n_{g_0}v(t)||_{L^{2}(M,g_0)}^{\frac{1}{2}}||\n_{g_0}v(t)||_{W^{1,2}(M,g_0)}^{\frac{1}{2}},
$$
and
$$
||\n_{g_0}v(t)||_{L^2(M,g_0)}\leq C||v(t)||^{\frac{1}{2}}_{L^2(M,g_0)}||v(t)||_{W^{2,2}(M,g_0)}^{\frac{1}{2}}.
$$
The last two formulas imply
\begin{equation}\label{eq:183v}
||\n_{g_0}v(t)||_{L^{\frac{8}{3}}(M,g_0)}\leq C||v(t)||_{L^2(M,g_0)}^{\frac{1}{4}}||v(t)||_{W^{2,2}(M,g_0)}^{\frac{3}{4}}.
\end{equation}
Hence \;$\eqref{eq:18u}$,\;$\eqref{eq:nunv}$\;and \;$\eqref{eq:183v}$\;imply
\begin{equation}\label{eq:nunv1}
\int_{M}(\n_{0}u(t),\n_{g_0}v(t))^2dV_{g_0}\leq C|||u(t)||_{W^{2,2}(M,g_0)}^{\frac{3}{2}}||u(t)||_{W^{4,2}(M,g_0)}^{\frac{1}{2}}||v(t)||_{L^2(M,g_0)}^{\frac{1}{2}}||v(t)||_{W^{2,2}(M,g_0)}^{\frac{3}{2}};
\end{equation}
Still from H\"{o}lder inequality we get
\begin{equation}\label{eq:luv}
 \int_{M}(\D_{g_0}u(t))^2v(t)^2dV_{g_0}\leq ||\D_{g_0}u(t)||_{L^4(M,g_0)}^2||v(t)||_{L^4(M,g_0)}^2.
\end{equation}
Furthermore using again interpolation inequality, we have
$$
||\D_{g_0}u(t)||_{L^{4}(M,g_0)}\leq C||\D_{g_0}||_{L^2(M,g_0)}^{\frac{1}{2}}||\D_{g_0}u(t)||_{W^{2,2}(M,g_0)}^{\frac{1}{2}}. 
$$
Thus, we obtain
\begin{equation}\label{eq:24u}
||\D_{g_0}u(t)||_{L^4(M,g_0)} \leq C||u(t)||_{W^{2,2}(M,g_0)}^{\frac{1}{2}}||u(t)||_{W^{4,2}(M,g_0)}^{\frac{1}{2}}.
\end{equation}
Hence \;$\eqref{eq:04v}$, \;$\eqref{eq:luv}$\;and \;$\eqref{eq:24u}$\;imply
\begin{equation}\label{eq:luv1}
 \int_{M}(\D_{g_0}u(t))^2v(t)^2dV_{g_0}\leq ||u(t)||_{W^{2,2}(M,g_0)}||u(t)||_{W^{4,2}(M,g_0)}||v(t)||_{L^2(M,g_0)}||v(t)||_{W^{2,2}(M,g_0)}.
\end{equation}
Now using the fact that \;$||u(t)||_{W^{2,2}(M,g_0)}\leq C$\;(see Lemma\;$\ref{eq:eh2b}$), we have that \;$\eqref{eq:lapar}$,\;$\eqref{eq:guv1}$ \;$\eqref{eq:nunv1}$\; and \;$\eqref{eq:luv1}$\;imply 
\begin{equation}\label{eq:esqp}
\begin{split}
 \int_{M}-4\D_{g_0}(v(t)e^{2u(t)})\D_{g_0}(v(t)e^{-2u(t)})dV_{g_0}\leq &-2\int_{M}(\D_{g_0}v(t))^2dV_{g_0}\\&+C\left(|u(t)||_{W^{4,2}(M,g_0)}||v(t)||_{L^2(M,g_0)}||v(t)||_{W^{2,2}(M,g_0)}\right)\\&+C\left(||u(t)||_{W^{4,2}(M,g_0)}^{\frac{1}{2}}||v(t)||_{L^2(M,g_0)}^{\frac{1}{2}}||v(t)||_{W^{2,2}(M,g_0)}^{\frac{3}{2}}\right)
\end{split}
\end{equation}
Next using the relation 
\begin{equation}\label{eq:iden}
\n_{g_0}(v(t)e^{2u(t)})\n_{g_0}(v(t)e^{-2u(t)}))=(|\n_{g_0}v(t)|^2-4|\n_{g_0}u(t)|^2v(t)^2),
\end{equation}
we get 
\begin{equation}
 \int_{M}\left|\n_{g_0}(v(t)e^{2u(t)})\n_{g_0}(v(t)^2e^{-2u(t)}))\right|dV_{g_0}\leq ||\n_{g_0}v(t)||^2_{L^2(M,g_0)}+4\int_{M}|\n_{g_0}u(t)|^2v(t)^2dV_{g_0}
\end{equation}
Now using again H\"{o}lder inequality, we have
$$
\int_{M}|\n_{g_0}u(t)|^2v(t)^2dV_{g_0}\leq C||\n_{g_0}u(t)||_{L^4(M,g_0)}^2||v(t)||_{L^4(M,g_0)}^2.
$$
Furthermore, applying the Sobolev embedding theorem, we get
$$
||\n_{g_0}u(t)||_{L^4(M,g_0)}\leq ||u(t)||_{W^{2,2}(M,g_0)}.
$$
Thus\;$\eqref{eq:04v}$\;and the fact that \;$||u(t)||_{W^{2,2}(M,g_0)}\leq C$\;imply
\begin{equation*}
 \int_{M}|\n_{g_0}u(t)|^2v(t)^2dV_{g_0}\leq C||v(t)||_{W^{2,2}(M,g_0)}||v(t)||_{L^2(M,g_0)}.
\end{equation*}
Hence we obtain 
\begin{equation}\label{eq:scric0}
\begin{split}
\left| \int_{M}-\frac{8}{3}R_{g_0}\n_{g_0}(ve^{2u})\n_{g_0}(ve^{-2u})+8Ric_{g_0}(\n_{g_0}(ve^{2u}),\n_{g_0}(ve^{-2u}))dV_{g_0}\right|\leq  C ||v(t)||_{L^2(M,g_0)}||v(t)||_{W^{2,2}(M,g_0)}\\+C||\n_{g_0}v(t)||^2_{L^2(M,g_0)}.
\end{split}
\end{equation}
On the other hand using again interpolation inequality, we get
$$
||\n_{g_0}v(t)||^2_{L^2(M,g_0)}\leq C||v(t)||_{L^2(M,g_0)}||v(t)||_{W^{2,2}(M,g_0)}.
$$
So we have
\begin{equation}\label{eq:scric}
\left| \int_{M}-\frac{8}{3}R_{g_0}\n_{g_0}(ve^{2u})\n_{g_0}(ve^{-2u})+8Ric_{g_0}(\n_{g_0}(ve^{2u}),\n_{g_0}(ve^{-2u}))dV_{g_0}\right|\leq  C ||v(t)||_{L^2(M,g_0)}||v(t)||_{W^{2,2}(M,g_0)}.
\end{equation}
Next by using the expression of \;$P^{4}_{g_0}$\;one can easily check that the folowing holds
\begin{equation}
|P^{4}_{g_0}(e^{4u(t)})|\leq Ce^{4u(t)}\left(|\n_{g_0}^4u(t)|+|\n_{g_0}^2u(t)|+|\n_{g_0}u(t)||\n_{g_0}^3u(t)|+|\n_{g_0}^2u(t)|^2+|\n_{g_0}u(t)|^2+|\n_{g_0}u(t)|^4\right).
\end{equation}
Thus, we obtain
$$
|e^{-4u(t)}P^{4}_{g_0}(e^{4u(t)})|\leq C\left(|\n_{g_0}^4u(t)|+|\n_{g_0}^2u(t)|+|\n_{g_0}u(t)||\n_{g_0}^3u(t)|+|\n_{g_0}^2u(t)|^2+|\n_{g_0}u(t)|^2+|\n_{g_0}u(t)|^4\right).
$$
Now, taking the square of both side, integrating both sides, after taking the square root and using H\"{o}lder inequality and Sobolev embedding we get 
\begin{equation}\label{eq:topro}
\begin{split}
 ||e^{-4u}P^{4}_{g_0}(e^{4u(t)})||_{L^2(M,g_0)}\leq C\left(||u(t)||_{W^{4,2}(M,g_0)}+||u(t)||_{W^{1,4}(M,g_0)}||u(t)||_{W^{3,4}(M,g_0)}+||u(t)||_{W^{2,4}(M,g_0)}^2\right)\\+C||u(t)||_{W^{1,8}(M,g_0)}^4.
\end{split}
\end{equation}On the other hand, using Sobolev embedding theorem, we have
$$
||u(t)||_{W^{1,4}(M,g_0)}\leq C||u(t)||_{W^{2,2}(M,g_0)}
$$
and
$$
||u(t)||_{W^{3,4}(M,g_0)}\leq C||u(t)||_{W^{4,2}(M,g_0)}.
$$
From the fact that \;$-C\leq \ov{u(t)}_{g_0}\leq C$, we get 
$$
||u(t)||_{W^{2,4}(M,g_0)}\leq C||\D_{g_0}u(t)||_{L^4(M,g_0)}.
$$
Hence, using \;$\eqref{eq:24u}$\;we obtain
\begin{equation}\label{eq:sauveur}
||u(t)||_{W^{2,4}(M,g_0)}\leq C||u(t)||_{W^{2,2}(M,g_0)}^{\frac{1}{2}}||u(t)||_{W^{4,2}(M,g_0)}^{\frac{1}{2}}
\end{equation}
Thus \;$\eqref{eq:18u}$\;and\;$\eqref{eq:sauveur}$\; imply that\;$\eqref{eq:topro}$\;becomes
\begin{equation}
\begin{split}
 ||e^{-4u(t)}P^{4}_{g_0}(e^{4u(t)})||_{L^2(M,g_0)}\leq C\left(||u(t)||_{W^{4,2}(M,g_0)}+||u(t)||_{W^{2,2}( M,g_0)}||u(t)||_{W^{4,2}(M,g_0)}\right)\\+C\left(||u(t)||_{H^2(M,g_0)}^3||u(t)||_{W^{4,2}(M,g_0)}\right).
\end{split}
\end{equation}
Thus using again the fact that \;$||u(t)||_{W^{2,2}(M,g_0)}\leq C$, we get
\begin{equation}\label{eq:useful}
||e^{-4u(t)}P^{4}_{g_0}(e^{4u(t)})||_{L^2(M,g_0)}\leq C||u(t)||_{W^{4,2}(M,g_0)}.
\end{equation}
Now from Moser-Trudinger inequality we infer that 
\begin{equation}\label{eq:apm}
 \sup\{||e^{2u(t)}||_{L^4(M,g_0)},||e^{-2u(t)}||_{L^2(M,g_0)}\}\leq C.
\end{equation}
Using Young's inequality we get 
$$
\int_{M}|P^4_{g_0}(e^{4u(t)})|ve^{-2u(t)}\leq C||e^{-4u(t)}P^{4}_{g_0}(e^{4u(t)})||_{L^2(M,g_0)}||v||_{L^4(M,g_0)}||e^{2u(t)}||_{L^4(M,g_0)}
$$
Hence \;$\eqref{eq:04v}$\;\;$\eqref{eq:useful}$\;and\;$\eqref{eq:apm}$\;imply
\begin{equation}\label{eq:part1}
 \int_{M}|P^4_{g_0}(e^{4u(t)})|ve^{-2u(t)}\leq C||u(t)||_{W^{4,2}(M,g_0)}||v(t)||_{H^2(M,g_0)}.
\end{equation}
On the other hand using the definition of \;$P^3_{g_0}$\;one can check easily that the following hold
$$
|P^3_{g_0}(e^{4u(t)})|\leq Ce^{4u(t)}\left(|\frac{\partial \D_{g_0}u(t)}{\partial n_{g_0}}|+|\n_{\hat g_0}^2u(t)|+|\frac{\partial (|\n_{g_0}u(t)|^2)}{\partial n_{g_0}}|+|\n_{\hat g_0}u(t)|\right).
$$
Hence, using the same trick as above, we obtain
\begin{equation}\label{eq:ep3}
\begin{split}
||e^{-4u(t)}P^3_{g_0}(e^{4u(t)})||_{L^2(\partial M,\hat g_0)}\leq C\left(||\frac{\partial \D_{g_0}u(t)}{\partial n_{g_0}}||_{L^2(\partial M,\hat g_0)}+||\n_{\hat g_0}^2u(t)||_{L^2(\partial M,\hat g_0)}+||\frac{\partial (|\n_{g_0}u(t)|^2)}{\partial n_{g_0}}||_{L^2(\partial M,\hat g_0)}\right)\\+C||\n_{\hat g_0}u(t)||_{L^2(\partial M,\hat g_0)}.
\end{split}
\end{equation}
Now using trace Sobolev embedding we get 
\begin{equation}\label{eq:tnlap}
||\frac{\partial \D_{g_0}u(t)}{\partial n_{g_0}}||_{W^{\frac{1}{2},2}(\partial M,\hat g_0)}\leq C||\D_{g_0}u||_{W^{2,2}(M,g_0)},
\end{equation}
\begin{equation}\label{eq:tngu2}
 ||\frac{\partial (|\n_{g_0}u(t)|^2)}{\partial n_{g_0}}||_{W^{\frac{1}{2},2}(\partial M,\hat g_0)}\leq C|||\n_{g_0}u(t)|^2||_{W^{2,2}(M,g_0)},
\end{equation}
and 
\begin{equation}\label{eq:n2n1}
 ||\n_{\hat {g_0}}^2u(t)||_{L^2(\partial M,\hat g_0)}+||\n_{\hat {g_0}}u(t)||_{L^2(\partial M,hat g_0)}\leq C||u(t)||_{W^{2+\frac{1}{2},2}(\partial M,\hat g_0)}\leq C||u(t)||_{W^{3,2}(M,g_0)}.
\end{equation}
Hence \;$\eqref{eq:tnlap}$-$\eqref{eq:n2n1}$\;imply that \;$\eqref{eq:ep3}$\;becomes
\begin{equation}\label{eq:86}
 ||e^{-4u(t)}P^3_{g_0}(e^{4u(t)})||_{L^2(\partial M)}\leq C\left(||\D_{g_0}u||_{W^{2,2}(M,g_0)}+|||\n_{g_0}u(t)|^2||_{w^{2,2}(M,g_0)}+||u(t)||_{W^{3,2}(M,g_0)}\right).
\end{equation}
Now using Sobolev embedding, we obtain
\begin{equation}
 ||\D_{g_0}u||_{W^{2,2}(M,g_0)}\leq ||u||_{W^{4,2}(M,g_0)},
\end{equation}
\begin{equation}
 ||u(t)||_{W^{3,2}(M,g_0)}\leq C|||u(t)||_{W^{4,2}(M,g_0)}.
\end{equation}
Furthermore using H\"{o}lder inequality, we derive
\begin{equation}
 |||\n_{g_0}u(t)|^2||_{W^{2,2}(M,g_0)}\leq C\left(||\n_{g_0}^3u(t)||_{L^4(M,g_0)}||\n_{g_0}u(t)||_{L^4(M,g_0)}+||\n_{g_0}^2u(t)||^2_{L^4(M,g_0)}\right);
\end{equation}
hence we obtain
\begin{equation}
 |||\n_{g_0}u(t)|^2||_{W^{2,2}(M,g_0)}\leq\left(||u(t)||_{W^{3,4}(M,g_0)}||u(t)||_{W^{1,4}(M,g_0)}+||u(t)||^2_{W^{2,4}(M,g_0)}\right).
\end{equation}
Now using again Sobolev embedding we get
\begin{equation}
 ||u(t)||_{W^{3,4}(M,g_0)}\leq C||u(t)||_{W^{4,2}(M,g_0)};
\end{equation}
and
\begin{equation}\label{eq:93}
||u(t)||_{W^{1,4}(M,g_0)}\leq C||u(t)||_{W^{2,2}(M,g_0)}
\end{equation}
Thus from the fact that \;$||u(t)||_{W^{2,2}(M,g_0)}\leq C$, \;$\eqref{eq:sauveur}$\;and  \;$\eqref{eq:86}$-$\eqref{eq:93}$\;we infer 
\begin{equation}\label{eq:94}
 ||e^{-4u(t)}P^3_{g_0}(e^{4u(t)})||_{L^2(\partial M,\hat g_0)}\leq C||u(t)||_{W^{4,2}(M,g_0)}.
\end{equation}
 On the other hand using H\"{o}lder inequality we obtain
\begin{equation}\label{eq:part2}
 \int_{\partial M}|P^3_{g_0}(e^{4u(t)})|ve^{-2u(t)}dV_{g_0}\leq ||e^{-4u(t)}P^3_{g_0}(e^{4u(t)})||_{L^2(\partial M)}||v(t)||_{L^4(\partial M)}||e^{2u(t)}||_{L^{2}(\partial M)}.
\end{equation}
Next using Moser-Trudinger inequality, trace Sobolev embedding and Sobolev embedding we get
\begin{equation}\label{eq:96}
 ||e^{2|u(t)|}||_{L^{4}(\partial M,\hat g_0)}\leq C,
\end{equation}
and
\begin{equation}\label{eq:97}
 ||v(t)||_{L^4(\partial M,\hat g_0)}\leq C||v(t)||_{W^{1,4}(M,g_0)}\leq C||v(t)||_{W^{2,2}(M,g_0)}.
\end{equation}
Hence we have \;$\eqref{eq:94}$, \;$\eqref{eq:96}$\;'and \;$\eqref{eq:97}$\;imply that \;$\eqref{eq:part2}$\;becomes
\begin{equation}\label{eq:part3}
 \int_{\partial M}|P^3_{g_0}(e^{4u(t)})|ve^{-2u(t)}dS_{g_0}\leq C||v(t)||_{W^{2,2}(M,g_0)}||u(t)||_{W^{4,2}(M,g_0)}.
\end{equation}
Now let us estimate
$$
\left|\int_{M}\left[ (4e^{4u(t)}\frac{\ov{Q_{g(t)}}}{\ov F}F-4Q_{g_0})P^{4,3}_{g_0}(v(t)e^{-2u(t)})\right]dV_{g_0}\right|
$$
Using the self-adjointness of \;$P^{4,3}_{g_0}$, we have 
\begin{equation}
\begin{split}
\int_{M}\left[ (4e^{4u(t)}\frac{\ov{Q_{g(t)}}}{\ov F}F-4Q_{g_0})P^{4,3}_{g_0}(v(t)e^{-2u(t)})\right]dV_{g_0}=\int_{M} 4\frac{\ov{Q_{g(t)}}}{\ov F}P^{4,3}_{g_0}(e^{4u(t)}F)(v(t)e^{-2u(t)})dV_{g_0}\\-\int_{M}4P^{4,3}_{g_0}(Q_{g_0})(v(t)e^{-2u(t)})dV_{g_0}
\end{split}
\end{equation}
Thus we obtain 
\begin{equation}
\begin{split}
\left|\int_{M}\left[ (4e^{4u(t)}\frac{\ov{Q_{g(t)}}}{\ov F}F-4Q_{g_0})P^{4,3}_{g_0}(v(t)e^{-2u(t)})\right]dV_{g_0}\right|\leq \left|\int_{M} 4\frac{\ov{Q_{g(t)}}}{\ov F}P^{4,3}_{g_0}(e^{4u(t)}F)(v(t)e^{-2u(t)})dV_{g_0}\right|\\+\left|\int_{M}4P^{4,3}_{g_0}(Q_{g_0})(v(t)e^{-2u(t)})dV_{g_0}\right|.
\end{split}
\end{equation}
Now, using H\"{o}lder inequality, we get
$$
\left|\int_{M}4P^{4,3}_{g_0}(Q_{g_0})(v(t)e^{-2u(t)})dV_{g_0}\right|\leq ||v(t)||_{L^4(M,g_0)}||e^{-2u(t)}||_{L^{4}(M,g_0)}.
$$
Hence using, we obtain
$$
\left|\int_{M}4P^{4,3}_{g_0}(Q_{g_0})(v(t)e^{-2u(t)})dV_{g_0}\right|\leq ||v(t)||_{L^2(M,g_0)}^{\frac{1}{2}}||v(t)||_{W^{2,2}(M,g_0)}^{\frac{1}{2}}.
$$
Next let us estimate\;$\left|\int_{M} P^{4,3}_{g_0}(e^{4u(t)}F)(v(t)e^{-2u(t)})dV_{g_0}\right|$.\\
First of all we have
\begin{equation}
 \begin{split}
\int_{M} P^{4,3}_{g_0}(e^{4u(t)}F)(v(t)e^{-2u(t)})dV_{g_0}=\int_{M} \left(P^{4,3}_{g_0}(e^{4u(t)}F)-FP^{4,3}_{g_0}(e^{4u(t)})\right)(v(t)e^{-2u(t)})dV_{g_0}\\+\int_{M} FP^{4,3}_{g_0}(e^{4u(t)})(v(t)e^{-2u(t)})dV_{g_0}.
\end{split}
\end{equation}
Thus, we get
\begin{equation}
\begin{split}
\left|\int_{M} P^{4,3}_{g_0}(e^{4u(t)}F)(v(t)e^{-2u(t)})dV_{g_0}\right|\leq \left|\int_{M} \left(P^{4,3}_{g_0}(e^{4u(t)}F)-FP^{4,3}_{g_0}(e^{4u(t)})\right)(v(t)e^{-2u(t)})dV_{g_0}\right|\\+\left|\int_{M} FP^{4,3}_{g_0}(e^{4u(t)})(v(t)e^{-2u(t)})dV_{g_0}.\right|
\end{split}
\end{equation}
Now using H\"{o}lder inequality. \;$\eqref{eq:part1}$, and \;$\eqref{eq:part3}$, we have that the second term in the left hand side of the above inequality can be estimated as follows
$$
\left|\int_{M} FP^{4,3}_{g_0}(e^{4u(t)})(v(t)e^{-2u(t)})dV_{g_0}.\right|\leq C||v(t)||_{W^{2,2}(M,g_0)}||u(t)||_{W^{4,2}(M,g_0)}.
$$
Next using again Holder inequality, we infer that the first term can be estimated in the following way
\begin{equation}
 \begin{split}
\left|\int_{M} \left(P^{4,3}_{g_0}(e^{4u(t)}F)-FP^{4,3}_{g_0}(e^{4u(t)})\right)(v(t)e^{-2u(t)})dV_{g_0}\right|\leq ||\left(P^{4,3}_{g_0}(e^{4u(t)}F)-FP^{4,3}_{g_0}(e^{4u(t)})\right)||_{L^{2}(M,g_0)}||v(t)||_{L^4(M,g_0)}\\||e^{-2u(t)}||_{L^4(M,g_0)}.
\end{split}
\end{equation}
On the other hand, using {\em commutators} formula in pseudodifferential calculus, we have
$$
||\left(P^{4,3}_{g_0}(e^{4u(t)}F)-FP^{4,3}_{g_0}(e^{4u(t)})\right)||_{L^{2}(M,g_0)}\leq C||e^{4u(t)}||_{W^{3,2}(M,g_0)}.
$$
moreover using the same argument as above one can check easily that the following holds
$$
||e^{4u(t)}||_{W^{3,2}(M,g_0)}\leq ||u(t)||_{W^{4,2}(M,g_0)}.
$$
Thus we derive 
\begin{equation}
\left|\int_{M} \left(P^{4,3}_{g_0}(e^{4u(t)}F)-FP^{4,3}_{g_0}(e^{4u(t)})\right)(v(t)e^{-2u(t)})dV_{g_0}\right|\leq ||u(t)||_{W^{4,2}(M,g_0)}||v(t)||_{L^2(M,g_0)}^{\frac{1}{2}}||v(t)||_{W^{2,2}(M,g_0)}^{\frac{1}{2}}.
\end{equation}
Hence combining all, and the fact that \;$\ov{Q_{g(t)}}$, is invariant with respect to \;$t$\;we get
\begin{equation}
\int_{M}\left[ (4e^{4u(t)}\frac{\ov{Q_{g(t)}}}{\ov F}F-4Q_{g_0})P^{4,3}_{g_0}(v(t)e^{-2u(t)})\right]dV_{g_0}\leq C||v(t)||_{W^{2,2}(M,g_0)}||u(t)||_{W^{4,2}(M,g_0)}.
\end{equation}
Now let us estimate $8\int_{\partial M}L_{g_0}(\n_{\hat g_0}(v(t)e^{2u(t)}),\n_{\hat g_0}(v(t)e^{-2u(t)}))dS_{g_0}$. Using the identity \;$\eqref{eq:iden}$\;(with \;$g$\;replaced by \;$\tilde g$) we get
$$
\left|8\int_{\partial M}L_{g_0}(\n_{\tilde g}(ve^{2u(t)}),\n_{\tilde g}(ve^{-2u(t)}))dS_{g_0}\right|\leq C\int_{\partial M}|\n_{\tilde g}v(t)|^2dS_{g_0}+C\int_{\partial M}|\n_{\tilde g}u(t)|^2v(t)^2.
$$
Now for \;$\epsilon>0$\;small enough, we have by Lemma 2.3 in \cite{cq2}
\begin{equation}
 \int_{M}|\n_{\tilde {g_0}}v(t)|^2dS_{g_0}\leq \epsilon \int_{M}(\D_{g_0}v(t))^2+C_{\epsilon}||v(t)||_{W^{1,2}}^2.
\end{equation}
Hence using Sobolev embedding we get
\begin{equation}
 \int_{M}|\n_{\tilde {g_0}}v(t)|^2dS_{g_0}\leq \epsilon \int_{M}(\D_{g_0}v(t))^2+C_{\epsilon}||v(t)||_{L^2(M,g_0)}||v(t)||_{W^{2,2}(M,g_0)}.
\end{equation}
On the other hand, using again H\"{o}lder 	
inequality we get
\begin{equation}
 \int_{\partial M}|\n_{\hat g_0}u(t)|^2v(t)^2\leq ||\n_{\hat g_0}u(t)||^2_{L^6(\partial M, \hat g_0)}||v(t)||^2_{L^3(\partial M,\hat g_0)}.
\end{equation}
Now using trace Sobolev embedding and sobolev emdedding we derive
\begin{equation}
 ||\n_{\tilde{g_0}}u(t)||^2_{L^6(\partial M,\hat g_0)}\leq C||\n_{\tilde {g_0}}u(t)||^2_{W^{1,2}(\partial M,\hat g_0)}\leq C||u(t)||_{W^{3,3}(M,g_0)}^2;
\end{equation}
On the other hand using interpolation, it is easily seen that the following holds
$$
||u(t)||_{W^{3,3}(M,g_0)}^2\leq C||u(t)||_{W^{2,2}(M,g_0)}^{\frac{1}{3}}||u(t)||^{\frac{5}{3}}_{W^{4,2}(M,g_0)}
$$
Hence we obtain
$$
||\n_{\tilde{g_0}}u(t)||^2_{L^6(\partial M,\hat g_0)}\leq C||u(t)||^{\frac{1}{3}}_{W^{2,2}(M,g_0)}||u(t)||^{\frac{5}{3}}_{W^{4,2}(M,g_0)}.
$$
Now using trace Sobolev embedding, we get 
\begin{equation}
 ||v(t)||^2_{L^3(\partial M,\hat g_0)}\leq C||v(t)||^2_{W^{1,3}(M,g_0)}.
\end{equation}
Next  using interpolation, we obtain
$$
||v(t)||^2_{W^{1,3}(M,g_0)}\leq C||v(t)||_{L^2(M,g_0)}^{\frac{1}{3}}||v(t)||_{W^{2,2}(M,g_0)}^{\frac{5}{3}}.
$$
Hence we arrive to
$$
 ||v(t)||^2_{L^3(\partial M,\hat g_0)}\leq C ||v(t)||_{L^2(M,g_0)}^{\frac{1}{3}}||v(t)||_{W^{2,2}(M,g_0)}^{\frac{5}{3}}.
$$
Thus using the fact that \;$||u(t)||_{W^{2,2}(M,g_0)}\leq C$, we get
$$
 \int_{\partial M}|\n_{\hat g_0}u(t)|^2v(t)^2\leq C||u(t)||^{\frac{5}{3}}_{W^{4,2}(M,g_0)}||v(t)||_{L^2(M,g_0)}^{\frac{1}{3}}||v(t)||_{W^{2,2}(M,g_0)}^{\frac{5}{3}}
$$
Hence collecting all, we get
\begin{equation}
\begin{split}
 \left|8\int_{\partial M}L_{g_0}(\n_{\tilde g}(v(t)e^{2u(t)}),\n_{\tilde g}(v(t)e^{-2u(t)}))dS_{g_0}\right|\leq \epsilon \int_{M}(\D_{g_0}v(t))^2+C_{\epsilon}||v(t)||_{L^2(M,g_0)}||v(t)||_{W^{2,2}(M,g_0)}\\+C||u(t)||^{\frac{5}{3}}_{W^{4,2}(M,g_0)}||v(t)||_{L^2(M,g_0)}^{\frac{1}{3}}||v(t)||_{W^{2,2}(M,g_0)}^{\frac{5}{3}}.
\end{split}
\end{equation}
Hence we arrive to 
\begin{equation*}
\begin{split}
\frac{d}{d t}\left(\int_{M}(P^{4,3}_{g_0}(u(t)))^2dV_{g_0}\right)\leq -2\int_{M}(\D_{g_0}v(t))^2+ C ||v(t)||_{W^{2,2}(M,g_0)}||v(t)||_{L^2(M,g_0)}
\\+C||v(t)||_{W^{2,2}(M,g_0)}||u(t)||_{W^{4,2}(M,g_0)}\\+C\left(||u(t)||_{W^{4,2}(M,g_0)}||v||_{L^2(M,g_0)}||v(t)||_{W^{2,2}(M,g_0)}+||u(t)||_{W^{4,2}(M,g_0)}^{\frac{1}{2}}||v(t)||_{L^2(M,g_0)}^{\frac{1}{2}}||v(t)||_{W^{2,2}(M,g_0)}^{\frac{3}{2}}\right)\\+ \epsilon \int_{M}(\D_{g_0}v(t))^2+C_{\epsilon}||v(t)||_{L^2(M,g_0)}||v(t)||_{W^{2,2}(M,g_0)}+C||u(t)||^{\frac{5}{3}}_{W^{4,2}(M,g_0)}||v(t)||^{\frac{1}{3}}_{L^2(M,g_0)}||v(t)||^{\frac{5}{3}}_{W^{2,2}(M,g_0)}.
\end{split}
\end{equation*}
From this and taking \;$\epsilon$\;so small, we have that by interpolation
\begin{equation}
\frac{d}{d t}\left(\int_{M}(P^{4,3}_{g_0}(u(t)))^2dV_{g_0}\right)\leq C(||v(t)||^2_{L^2(M,g_0)}+1)(||u(t)||^2_{W^{4,2}(M,g_0)}+1).
\end{equation}
Now using the fact that \;$P^{4,3}_{g_0}$\;is non-negative with trivial kernel  and the the relation \;$|\ov{u(t)}_{g_0}|\leq C$\;we derive the following
\begin{equation}
\int_{M}(P^{4,3}_{g_0}u(t))^2dV_0\geq C||u(t)||^2_{W^{4,2}}(M,g_0)-C.
\end{equation}
Furthermore, from the definition  of \;$v(t)$\; one can check easily that the following holds
\begin{equation}
||v(t)||^2_{L^2(M,g_0)}=\int_{M}(Q_{g(t)}-\frac{\ov{ Q_{g(t)}}}{\ov F}F)^2dV_{g(t)}.
\end{equation}
Hence we have that \;$\int_{M}(P^{4,3}_{g_0}u(t))^2dV_{g_0}+1$\;  satisfies the following differential inequality
\begin{equation}
\frac{d}{d t}\left(\int_{M}(P^{4,3}_{g_0}(u(t)))^2dV_{g_0}+1\right)\leq C(\int_{M}(Q_{g(t)}-\frac{\ov{ Q_{g(t)}}}{\ov F}F)^2dV_g+1)(\int_{M}(P^{4,3}_{g_0}(u(t)))^2dV_{g_0}+1).
\end{equation}
On the other hand, using the fact that \;$II_{Q,F}(u(t))$\;is bounded from below and \;$\frac{d II_{Q,F}(u(t))}{dt}=-4\int_{M}(Q_{g(t)}-\frac{\ov{ Q_{g(t)}}}{\ov F}F)^2dV_{g(t)}$\; we have that the following holds
\begin{equation}
\int_{0}^T\int_{M}(Q_{g(t)}-\ov{ Q_{g(t)}})^2dV_{g(t)}dt\leq C.
\end{equation}
Thus from the differential inequality that verifies \;$ \left(\int_{M}(P^{4,3}_{g_0}(u(t)))^2dV_{g_0}+1\right)$, we infer that
\begin{equation}
\int_{M}(P^{4,3}_{g_0}(u(t)))^2dV_{g_0}\leq C\;\;\;\text{for}\;\;\;0\leq t\leq T.
\end{equation}
So we obtain
\begin{equation*}
||u(t)-\ov{u(t)}||_{W^{4,2}(M,g_0)}\leq C \;\;\;\text{for}\;\;\;0\leq t< T.
\end{equation*}
Therefore, recalling that \;$-C\leq \ov{u(t)}_{g_0}\leq C$\;we arrive to
$$
||u(t)||_{W^{4,2}(M,g_0)}\leq C \;\;\;\text{for}\;\;\;0\leq t< T
$$
as desired. Hence the proof of the Lemma is concluded.
\end{pf}

\begin{pro}\label{eq:highcon}
Suppose that the assumptions of Theorem\;$\ref{eq:tqf}$\;holds and 
let \;$u$\;be the unique short-time solution to initial boundary value problem corresponding to \;$\eqref{eq:evu}$. Then for every \;$T>0$\;such that \;
$u$\;is defined on $[0, T[$, and and for every positive (relative) integer\;$k\geq 3$ there exists \;$C=C(T,k)>0$\;such that
\begin{equation}
||u(t)||_{W^{2k,2}(M)}\leq C\;\;\;\;\forall t\in [0, T[.
\end{equation}
\end{pro}
\begin{pf}
 We first derive a differential inequality for \;$\int_{M}\left((P^{4,3}_{g_0})^{\frac{k}{2}}(u(t))\right)^2dV_{g_0}$. To do so we compute \;$\frac{d}{dt}\left(\int_{M}((P^{4,3}_{g_0})^{\frac{k}{2}}(u(t)))^2dV_{g_0}\right)$. Using the rule of differentiation under the sign integral, we get
$$
\frac{d}{dt}\left(\int_{M}((P^{4,3}_{g_0})^{\frac{k}{2}}(u(t)))^2dV_{g_0}\right)=2\int_{M}(P^{4,3}_{g_0})^{\frac{k}{2}}(u(t))(P^{4,3}_{g_0})^{\frac{k}{2}}(\frac{\partial u(t)}{\partial t})dV_{g_0}.
$$
From the self-adjointness of \;$P^{4,3}_{g_0}$, we get
$$
\frac{d}{dt}\left(\int_{M}((P^{4,3}_{g_0})^{\frac{k}{2}}(u(t)))^2dV_{g_0}\right)=2\int_{M}(P^{4,3}_{g_0})^k(u(t))(\frac{\partial u(t)}{\partial t})dV_{g_0}.
$$
Next, using the eveolution equation for \;$u(t)$, we obtain
$$
\frac{d}{dt}\left(\int_{M}((P^{4,3}_{g_0})^{\frac{k}{2}}(u(t)))^2dV_{g_0}\right)=-\int_{M}(P^{4,3}_{g_0})^ku(t)\left(e^{-4u(t)}(P^4_{g_0}u(t)+2Q_{g_0})-2\frac{\ov{Q_{g(t)}}}{\ov F}F\right)dV_{g_0}.
$$
Since \;$\ov{Q_{g(t)}}$\;and \;$Vol_{g(t)}(M)$\;are invariant with respect to \;$t$, then we have
\begin{equation}
 \begin{split}
\frac{d}{dt}\left(\int_{M}((P^{4,3}_{g_0})^{\frac{k}{2}}(u(t)))^2dV_{g_0}\right)=-\int_{M}(P^{4,3}_{g_0})^ku(t)\left(e^{-4u(t)}(P^4_{g_0}u(t)+2Q_{g_0})\right)dV_{g_0}\\+O\left(||(P^{4,3}_{g_0})^{\frac{k}{2}}(u(t))||_{L^{2}(M,g_0)}\right).
\end{split}
\end{equation}
Expanding the left hand side of the later inequality, we obtain
\begin{equation}
 \begin{split}
\frac{d}{dt}\left(\int_{M}((P^{4,3}_{g_0})^{\frac{k}{2}}(u(t)))^2dV_{g_0}\right)=-\int_{M}(P^{4,3}_{g_0})^ku(t)e^{-4u(t)}P^4_{g_0}(u(t))dV_{g_0}-2\int_{M}(P^{4,3}_{g_0})^k(u(t))e^{-4u(t)}Q_{g_0}dV_{g_0}\\+O\left(||(P^{4,3}_{g_0})^{\frac{k}{2}}(u(t))||_{L^{2}(M,g_0)}\right).
\end{split}
\end{equation}
Using the fact that \;$P^3_{g_0}(u(t))=0$, we get
\begin{equation}
 \begin{split}
\frac{d}{dt}\left(\int_{M}((P^{4,3}_{g_0})^{\frac{k}{2}}(u(t)))^2dV_{g_0}\right)=-\int_{M}(P^{4,3}_{g_0})^ku(t)e^{-4u(t)}P^{4,3}_{g_0}(u(t))dV_{g_0}-2\int_{M}(P^{4,3}_{g_0})^k(u(t))e^{-4u(t)}Q_{g_0}dV_{g_0}\\+O\left(||(P^{4,3}_{g_0})^{\frac{k}{2}}(u(t))||_{L^{2}(M,g_0)}\right).
\end{split}
\end{equation}
From the self-adjointness of \;$P^{4,3}_{g_0}$, we have
\begin{equation}
 \begin{split}
\frac{d}{dt}\left(\int_{M}((P^{4,3}_{g_0})^{\frac{k}{2}}(u(t)))^2dV_{g_0}\right)=-\int_{M}(P^{4,3}_{g_0})^{\frac{k}{2}+\frac{1}{2}}(u(t))(P^{4,3}_{g_0})^{\frac{k}{2}-\frac{1}{2}}(e^{-4u(t)}P^{4,3}_{g_0}(u(t)))dV_{g_0}\\-2\int_{M}(P^{4,3}_{g_0})^{\frac{k}{2}}(u(t))(P^{4,3}_{g_0})^{\frac{k}{2}}(e^{-4u(t)}Q_{g_0})dV_{g_0}\\+O\left(||(P^{4,3}_{g_0})^{\frac{k}{2}}(u(t))||_{L^{2}(M,g_0)}\right).
\end{split}
\end{equation}
The latter expression can be rewritten in the following way
\begin{equation}
 \begin{split}
\frac{d}{dt}\left(\int_{M}((P^{4,3}_{g_0})^{\frac{k}{2}}(u(t)))^2dV_{g_0}\right)=-\int_{M}e^{-4u(t)}\left((P^{4,3}_{g_0})^{\frac{k}{2}+\frac{1}{2}}(u(t))\right)^2dV_{g_0}+O\left(||(P^{4,3}_{g_0})^{\frac{k}{2}}(u(t))||_{L^{2}(M,g_0)}\right)\\-\int_{M}(P^{4,3}_{g_0})^{\frac{k}{2}+\frac{1}{2}}(u(t))\left[(P^{4,3}_{g_0})^{\frac{k}{2}-\frac{1}{2}}\left(e^{-4(u(t)}P^{4,3}_{g_0}(u(t))\right)-e^{-4u(t)}(P^{4,3}_{g_0})^{\frac{k}{2}-\frac{1}{2}}((P^{4,3}_{g_0})u(t))\right]\\-2\int_{M}(P^{4,3}_{g_0})^{\frac{k}{2}}(u(t))(P^{4,3}_{g_0})^{\frac{k}{2}}(e^{-4u(t)}Q_{g_0})dV_{g_0}.\end{split}
\end{equation}
Using the same trick as above, we obtain
\begin{equation}
 \begin{split}
\frac{d}{dt}\left(\int_{M}((P^{4,3}_{g_0})^{\frac{k}{2}}(u(t)))^2dV_{g_0}\right)=-\int_{M}e^{-4u(t)}\left((P^{4,3}_{g_0})^{\frac{k}{2}+\frac{1}{2}}(u(t))\right)^2dV_{g_0}+O\left(||(P^{4,3}_{g_0})^{\frac{k}{2}}(u(t))||_{L^{2}(M,g_0)}\right)\\-\int_{M}(P^{4,3}_{g_0})^{\frac{k}{2}+\frac{1}{2}}(u(t))\left[(P^{4,3}_{g_0})^{\frac{k}{2}-\frac{1}{2}}\left(e^{-4(u(t)}P^{4,3}_{g_0}(u(t))\right)-e^{-4u(t)}(P^{4,3}_{g_0})^{\frac{k}{2}-\frac{1}{2}}((P^{4,3}_{g_0})u(t))\right]dV_{g_0}\\-2\int_{M}(P^{4,3}_{g_0})^{\frac{k}{2}}(u(t))\left[(P^{4,3}_{g_0})^{\frac{k}{2}}\left(e^{-4(u(t)}Q_{g_0}\right)-e^{-4u(t)}(P^{4,3}_{g_0})^{\frac{k}{2}}(Q_{g_0})\right]dV_{g_0}\\-2\int_{M}e^{-4u(t)}(P^{4,3}_{g_0})^{\frac{k}{2}}(u(t))(P^{4,3}_{g_0})^{\frac{k}{2}}(Q_{g_0})dV_{g_0}.\end{split}
\end{equation}
Now from the \;$W^{4,2}$-estimates, we derive from Sobolev embedding that 
$$
||u(t)||_{C^0(M)}\leq C.
$$
Thus using H\"{o}lder inequality, we have that the last term can be estimated as follows
$$
\left|2\int_{M}e^{-4u(t)}(P^{4,3}_{g_0})^{\frac{k}{2}}(u(t))(P^{4,3}_{g_0})^{\frac{k}{2}}(Q_{g_0})dV_{g_0}\right|\leq C||(P^{4,3}_{g_0})^{\frac{k}{2}}(u(t))||_{L^2(M,g_0)}.
$$
Still using H\"{o}lder inequality, we have that the following estimate holds for the third term
\begin{equation}
\begin{split}
\left|\int_{M}(P^{4,3}_{g_0})^{\frac{k}{2}+\frac{1}{2}}(u(t))\left[(P^{4,3}_{g_0})^{\frac{k}{2}-\frac{1}{2}}\left(e^{-4(u(t)}P^{4,3}_{g_0}(u(t))\right)-e^{-4u(t)}(P^{4,3}_{g_0})^{\frac{k}{2}-\frac{1}{2}}((P^{4,3}_{g_0})u(t))\right]dV_{g_0}\right|\leq \\C||(P^{4,3}_{g_0})^{\frac{k}{2}+\frac{1}{2}}(u(t))||_{L^{2}(M,g_0)}||\left[(P^{4,3}_{g_0})^{\frac{k}{2}-\frac{1}{2}}\left(e^{-4(u(t)}P^{4,3}_{g_0}(u(t))\right)-e^{-4u(t)}(P^{4,3}_{g_0})^{\frac{k}{2}-\frac{1}{2}}((P^{4,3}_{g_0})u(t))\right]||_{L^{2}(M,g_0)}
\end{split}
\end{equation}
On the other hand. still using H\"{o}lder inequality, we have the following estimate for the fourth term
\begin{equation}
\begin{split}
\left|2\int_{M}(P^{4,3}_{g_0})^{\frac{k}{2}}(u(t))\left[(P^{4,3}_{g_0})^{\frac{k}{2}}\left(e^{-4(u(t)}Q_{g_0}\right)-e^{-4u(t)}(P^{4,3}_{g_0})^{\frac{k}{2}}(Q_{g_0})\right]dV_{g_0}\right|\leq \\C||(P^{4,3}_{g_0})^{\frac{k}{2}}(u(t))||_{L^{2}(M,g_0)}||\left[(P^{4,3}_{g_0})^{\frac{k}{2}}\left(e^{-4(u(t)}Q_{g_0}\right)-e^{-4u(t)}(P^{4,3}_{g_0})^{\frac{k}{2}}(Q_{g_0}))\right]||_{L^{2}(M,g_0)}
\end{split}
\end{equation}
Now using {\em commutator} formulas in the theory of pseudodifferential operators, we get
\begin{equation}
 \begin{split}
||\left[(P^{4,3}_{g_0})^{\frac{k}{2}}\left(e^{-4(u(t)}Q_{g_0}\right)-e^{-4u(t)}(P^{4,3}_{g_0})^{\frac{k}{2}}(Q_{g_0}))\right]||_{L^{2}(M,g_0)}\leq \\C\sum_{k_1\geq 1,k_2\geq 0,k_1+k_2\leq 2k}||u(t)||_{W^{k_1,2}(M,g_0)}||Q_{g_0}||_{W^{k_2,2}(M,g_0)} \leq C||(P^{4,3}_{g_0})^{\frac{k}{2}}(u(t))||_{L^2(M,g_0)}.
 \end{split}
\end{equation}
and
\begin{equation}
 \begin{split}
 ||\left[(P^{4,3}_{g_0})^{\frac{k}{2}-\frac{1}{2}}\left(e^{-4(u(t)}P^{4,3}_{g_0}(u(t))\right)-e^{-4u(t)}(P^{4,3}_{g_0})^{\frac{k}{2}-\frac{1}{2}}((P^{4,3}_{g_0})u(t))\right]||_{L^{2}(M,g_0)}\leq \\
C\sum_{k_1\geq 1,k_2\geq 4,k_1+k_2\leq2k+2}||u(t)||_{W^{k_1,2}(M,g_0)}||u(t)||_{W^{k_2,2}(M,g_0)}.
 \end{split}
\end{equation}
Nex using interpolation, we obtain
$$
||u(t)||_{W^{k_1,2}(M,g_0)}\leq C||u(t)||_{L^2(M,g_0)}^{\frac{2k+2-k_1}{2k+2}}||u(t)||_{W^{2k+2,2}(M,g_0)}^{\frac{k_1}{2k+2}}
$$
and
$$
||u(t)||_{W^{k_2,2}(M,g_0)}\leq C||u(t)||_{W^2(M,g_0)}^{\frac{2k+2-k_2}{2k}}||u(t)||_{W^{2k+2,2}(M,g_0)}^{\frac{k_2-2}{2k}}
$$
Thus using again the fact that \;$||u(t)||_{W^{2,2}(M)}\leq C$, we obtain
\begin{equation}
 \begin{split}
 ||\left[(P^{4,3}_{g_0})^{\frac{k}{2}-\frac{1}{2}}\left(e^{-4(u(t)}P^{4,3}_{g_0}(u(t))\right)-e^{-4u(t)}(P^{4,3}_{g_0})^{\frac{k}{2}-\frac{1}{2}}((P^{4,3}_{g_0})u(t))\right]||_{L^{2}(M,g_0)}\leq \\
C\sum_{k_1\geq 1,k_2\geq 4,k_1+k_2\leq2k+2}||u(t)||_{W^{2k+2,2}(M,g_0)}^{\frac{k_1}{2k+2}+\frac{k_2-2}{2k}}
 \end{split}
\end{equation}
Collecting all, we get
\begin{equation}
\begin{split}
 \frac{d}{dt}\left(\int_{M}((P^{4,3}_{g_0})^{\frac{k}{2}}(u(t)))^2dV_{g_0}\right)\leq-\frac{1}{C}||u(t)||^{2}_{W^{2k+2}(M,g_0)}C\left(||u(t)||_{W^{2k,2}(M,g_0)}+||u(t)||^2_{W^{2k,2}(M,g_0)}+1\right)\\+C\sum_{k_1\geq 1,k_2\geq 4,k_1+k_2\leq2k+2}||u(t)||_{W^{2k+2,2}(M,g_0)}^{1+\frac{k_1}{2k+2}+\frac{k_2-2}{2k}}
\end{split}
\end{equation}
Thus, we have that \;$\int_{M}\left((P^{4,3}_{g_0})^{\frac{k}{2}}(u(t)\right)^2dV_{g_0}$\;verifies the following differential inequality
$$
\frac{d}{dt}\left[\int_{M}\left((P^{4,3}_{g_0})^{\frac{k}{2}}(u(t)\right)^2dV_{g_0}\right]\leq C\left(\int_{M}\left((P^{4,3}_{g_0})^{\frac{k}{2}}(u(t)\right)^2dV_{g_0}+1\right).
$$
Hence we obtain
$$
||u(t)||_{W^{2k,2}(M,g_0)}\leq C
$$
as desired.
\end{pf}

\subsection{Global existence and exponential convergence}
This subsection deals with the global exitence and the  convergence of the flow. Firts of all, the Proposition\;$\ref{eq:highcon}$\; rule out the formation of singularities in finite time. Hence it ensure the long-time existence. Now it remains to prove the convergence. To do so we start by giving the following Lemma.
\begin{lem}
Under the assumptions of Theorem\;$\ref{eq:tqf}$, seeting
\begin{equation*}
x(t)=\int_{M}(Q_{g(t)}-\frac{\ov{Q_{g(t)}}}{\ov F}F)^2dV_{g(t)},
\end{equation*}
we have that 
\begin{equation*}
x(t)\rightarrow 0\;\;\text{as}\;\;t\rightarrow +\infty.
\end{equation*}
\end{lem}
\begin{pf}
First of all using the fact that\;$II_{Q,F}$\;is bounded below and Proposition\;$\ref{eq:eh2b}$\; we have that \;$x(t)$\;is integrable. Hence  for \;$\epsilon>0$,\;$\delta>0$\;fixed, there exist a time \;$t_0$\;such that 
$$x(t_0)\leq \epsilon .$$
and
$$\int_{t_0}^{+\infty}x(t)dt\leq \delta.$$
We claim that $$x(t)\leq 3\epsilon\;\;\forall t\geq t_0.$$
To prove the claim we argue by contradiction. So let us suppose that the claim does not hold and let us argue for a contradiction. We define$$t_1=\inf\{t\geq t_0:\;\;x(t)\geq3\epsilon\}.$$
Therefore by definition \;$t_1$, we have that the following holds
\begin{equation}
x(t)\leq 3\epsilon\;\;\forall t\in [t_0, t_1].
\end{equation}
Hence from the fact that \;$\int_{M}Q_{g(t)}dV_{g(t)}$\;is invariant and also\;$\ov{Q_{g(t)}}$\;(with respect to \;$t$), we infer that 
\begin{equation}
\int_{M}Q_{g(t)}^2dV_{g(t)}\leq C \;\;\forall t\in[t_0, t_1].
\end{equation}
Therefore from the equation of transformation of \;$Q$-curvature, we derive
\begin{equation}
\int_{M}e^{-8u(t)}(P^4_{g_0}u(t)+2Q_{g_0})^2dV_{g_0}\leq C\;\;\;\forall t\in [t_0,t_1].
\end{equation}
On the other hand from the Moser-Trudinger inequality we have
\begin{equation}
\int_{M}e^{24u(t)}dV_{g_0}\leq C\;\;\;\forall \;\;t\in [t_0,t_1].
\end{equation}
Moreover using H\"older inequality we get
\begin{equation}
\int_{M}|P^4_{g_0}u(t)+2Q_{g_0}|^{\frac{3}{2}}dV_{g_0}\leq\left( \int_{M}e^{-8u(t)}(P^4_{g_0}u(t)+2Q_{g_0})dV_{g_0}\right)^{\frac{3}{4}}\left(\int_{M}e^{24u(t)}dV_{g_0}\right)^{\frac{1}{4}}.
\end{equation}
Therefore we obtain
\begin{equation}
\int_{M}|P^4_{g_0}u(t)|^{\frac{3}{2}}dV_{g_0}\leq C\;\;\;\forall \;t\in [t_0,t_1].
\end{equation}
So using the regularity result in Lemma\;$\ref{eq:reg}$, we get
\begin{equation}
||u(t)||_{W^{4,\frac{3}{2}}(M,g_0)}\leq C \;\;\forall \;t\in [t_0,t_1].
\end{equation}
Hence from Sobolev enbedding theroem we infer that
\begin{equation}
||u(t)||_{C^0(M)}\leq C \;\;\forall \;t\in [t_0,t_1].
\end{equation}
We remark here that \;$C$\;does not depend on \;$\epsilon,\delta,t_0,t_1$.\\
Now let us compute\;$\frac{d}{dt}\left(\int_{M}(Q_{g(t)}-\frac{\ov{Q_{g(t)}}}{\ov F}F)^2dV_{g(t)}\right)$. We have by differentiation rule under the sign integral, the evolution of the volume form , of the \;$Q$-curvature and the conformal factor that there holds
\begin{equation}
\begin{split}
\frac{d}{dt}\left(\int_{M}(Q_{g(t)}-\frac{\ov{Q_{g(t)}}}{\ov F}F)^2dV_{g(t)}\right)=8\int_{M}Q_{g(t)}(Q_{g(t)}-\frac{\ov{Q_{g(t)}}}{\ov F}F)^2dV_g-4\int_{M}(Q_{g(t)}-\frac{\ov{Q_{g(t)}}}{\ov F}F)^3dV_{g(t)}\\-\int_{M}(Q_{g(t)}-\frac{\ov{Q_{g(t)}}}{\ov F}F)P^4_{g(t)}(Q_{g(t)}-\frac{\ov{Q_{g(t)}}}{\ov F}F)dV_{g(t)}-8\frac{\ov{Q_{g(t)}}}{Vol_{g(t)}(M)}\left(\int_{M}\frac{F}{\ov F}(Q_{g(t)}-\frac{\ov{Q_{g(t)}}}{\ov F}F)dV_{g(t)}\right)^2.
\end{split}
\end{equation}
Hence using the expression of \;$P^{4,3}_{g(t)}$\;we arrive to 
\begin{equation}
\begin{split}
\frac{d}{dt}\left(\int_{M}(Q_{g(t)}-\frac{\ov{Q_{g(t)}}}{\ov F}F)^2dV_{g(t)}\right)=8\int_{M}Q_{g(t)}(Q_{g(t)}-\frac{\ov{Q_{g(t)}}}{\ov F}F)^2dV_g-4\int_{M}(Q_{g(t)}-\frac{\ov{Q_{g(t)}}}{\ov F}F)^3dV_{g(t)}\\-\int_{M}(Q_{g(t)}-\frac{\ov{Q_{g(t)}}}{\ov F}F)P^{4,3}_{g(t)}(Q_{g(t)}-\frac{\ov{Q_{g(t)}}}{\ov F}F)dV_{g(t)}+2\int_{\partial M}(Q_{g(t)}-\frac{\ov{Q_{g(t)}}}{\ov F}F)P^{3}_{g(t)}(Q_{g(t)}-\frac{\ov{Q_{g(t)}}}{\ov F}F)dS_{g(t)}\\-8\frac{\ov{Q_{g(t)}}}{Vol_{g(t)}(M)}\left(\int_{M}\frac{F}{\ov F}(Q_{g(t)}-\frac{\ov{Q_{g(t)}}}{\ov F}F)dV_{g(t)}\right)^2.
\end{split}
\end{equation}
On the other hand, using Gagliardo-Nirenberg inequality, we obtain 
\begin{equation}
\int_{M}|Q_{g(t)}-\frac{\ov{Q_{g(t)}}}{\ov F}F|^3dV_{g_0}\leq C\left(\int_{M}(Q_{g(t)}-\frac{\ov{Q_{g(t)}}}{\ov F}F)^2dV_{g_0}\right)||Q_{g(t)}-\frac{\ov{Q_{g(t)}}}{\ov F}F||_{W^{2,2}(M, g_0)}.
\end{equation}
So from Lemma\;$\ref{eq:neq}$, we infer that
\begin{equation}
\int_{M}(Q_{g(t)}-\frac{\ov{Q_{g(t)}}}{\ov F}F)^3dV_{g_0}\leq C\int_{M}(Q_{g(t)}-\frac{\ov{Q_{g(t)}}}{\ov F}F)^2dV_{g_0}\left(\int_{M}(Q_{g(t)}-\frac{\ov{Q_{g(t)}}}{\ov F}F)P^{4,3}_{g_0}(Q_{g(t)}-\frac{\ov{Q_{g(t)}}}{\ov F}F)dV_{g_{0}}\right)^{\frac{1}{2}}.
\end{equation}
Using  the fact that \;$||u(t)||_{C^0(M)}\leq C$\;and the conformal invariance of the Paneitz operator and the Chang-Qing one, we get 
\begin{equation}
\int_{M}|Q_{g(t)}-\frac{\ov{Q_{g(t)}}}{\ov F}F|^3dV_{g(t)}\leq C\int_{M}(Q_{g(t)}-\frac{\ov{Q_{g(t)}}}{\ov F}F)^2dV_{g(t)}\left(\int_{M}(Q_{g(t)}-\frac{\ov{Q_{g(t)}}}{\ov F}F)P^{4,3}_{g(t)}(Q_{g(t)}-\frac{\ov{Q_{g(t)}}}{\ov F}F)dV_{g(t)}\right)^{\frac{1}{2}}
\end{equation}
Now let \;$\gamma>0$\;small, using \;$\gamma$-Cauchy inequality we obtain
\begin{equation}
 \begin{split}
\int_{M}|Q_{g(t)}-\frac{\ov{Q_{g(t)}}}{\ov F}F|^3dV_{g(t)}\leq \gamma\left(\int_{M}(Q_{g(t)}-\frac{\ov{Q_{g(t)}}}{\ov F}F)P^{4,3}_{g(t)}(Q_{g(t)}-\frac{\ov{Q_{g(t)}}}{\ov F}F)dV_{g(t)}\right)\\+C_{\gamma}\left(\int_{M}(Q_{g(t)}-\frac{\ov{Q_{g(t)}}}{\ov F}F)^2dV_{g(t)}\right)^2.
\end{split}
\end{equation}

On the other hand , we have also
\begin{equation*}
\begin{split}
2\int_{\partial M}(Q_{g(t)}-\frac{\ov{Q_{g(t)}}}{\ov F}F)P^{3}_{g(t)}(Q_{g(t)}-\frac{\ov{Q_{g(t)}}}{\ov F}F)\leq  \gamma\left(\int_{M}(Q_{g(t)}-\frac{\ov{Q_{g(t)}}}{\ov F}F)P^{4,3}_{g(t)}(Q_{g(t)}-\frac{\ov{Q_{g(t)}}}{\ov F}F)dV_{g(t)}\right)\\+C_{\gamma}\left(\int_{M}(Q_{g(t)}-\frac{\ov{Q_{g(t)}}}{\ov F}F)^2dV_{g(t)}\right)^2.
\end{split}
\end{equation*}
Hence recollecting all, we obtain that \;$x(t)$\;satisfies the following differential inequality
\begin{equation}\label{eq:toint}
\frac{d}{dt}\left(x(t)\right)\leq Cx(t)^2+Cx(t).
\end{equation}
We remark that here \;$C$ depends on \;$\gamma$\;but not on \;$\epsilon$\;and \;$\delta$.\\
Now integrating both sides of \;$\eqref{eq:toint}$, we get
\begin{equation}
2\epsilon\leq x(t_1)-x(t_0)\leq C\int_{t_0}^{t_1}(x(t)^2+x(t))dt\leq C(1+\epsilon)\delta.
\end{equation}
Recalling that \;$C$\;does not depend on \;$\epsilon$\;and \;$\delta$, we reach a contradiction.
Therefore we have
\begin{equation}
\lim_{t\rightarrow+\infty}x(t)=0.
\end{equation}
as desired.
\end{pf}
\begin{pro}
Under the assumptions of Theorem\;$\ref{eq:tqf}$, the {\em eternal} solution to the initial boundary value problem corresponding to\;$\eqref{eq:qflowbf}$\; with initial data \;$g_0$\;is uniformly bounded in every \;$C^k$\;and converges to a smooth metric \;$g_{\infty}=e^{2u_{\infty}}$\;(conformal to \;$g_0$). Moreover the \;$Q$-curvature \;$Q_{\infty}$\;of the limiting metric \;$g_{\infty}$\;satisfy\;$Q_{\infty}=\frac{\ov {Q_{\infty}}}{\ov F}F$, and its \;$T$-curvature and mean curvature vanish.
\end{pro}
\begin{pf}
Arguing as in the Lemma above, we derive the following two facts
$$
\int_{M}e^{-8u(t)}(2Q_{g_0}+P^4_{g_0}u(t))^2dV_{g_0}\leq C.
$$
$$
||u(t)||_{C^0}\leq C\;.
$$
for every non-negative \;$t$. Hence we obtain
$$
\int_{M}(2Q_{g_0}+P^4_{g_0}u(t))^2dV_{g_0}\leq C,
$$
which implies
$$
\int_{M}(P^4_{g_0}u(t))^2dV_{g_0}\leq C
$$
Thus using the regularity result in Lemma\;$\ref{eq:reg}$, we derive 
\begin{equation}
||u(t)||_{W^{4,2}(M,g_0)}\leq C\;\;\forall t\geq 0.
\end{equation}
From this last estimate and arguing as in Proposition\;$\ref{eq:highcon}$, we obtain
\begin{equation}
||u(t)||_{w^{k,2}(M,g_0)}\leq C\;\;\forall t\geq 0\;\;\forall k\geq 4.
\end{equation}
Hence the boundedness in every \;$C^k$\;of \;$g(t)$\;follows.\\
Now, we remark the flow is the gradient flow of the  geometric functional \;$II_{Q,F}$. Hence from the fact that \;$II_{Q,F}$\;is real analytic, we have that the convergence follows directly from Leon Simon result, see \cite{si}.
\end{pf}
\section{Proof of Theorem\;$\ref{eq:ttf}$}
As already remarked at the Introduction, in this small Section, we give the proof of the short-time existence for the initial boundary value problem corresponding to \;$\eqref{eq:tflowbs}$.\\

First of all, if we write the evolving metric metric as \;$g(t)=e^{2u(t)}g_0$, we have that \;$u(t)$\;satifies the following evolution equation on\;$\partial M$
$$
\frac{\partial u(t)}{\partial t}=-\left(e^{-3u(t)}(P^3_{g_0}(u(t))+T_{g_0})-\frac{\ov T_{g(t)}}{\ov S}S\right).
$$
Moreover, we have also
$$
P^4_{g_0}(u(t))=0,\;\;\text{on}\;\;M\;\;\;\;\frac{\partial u(t)}{\partial n_{g_0}}=0,\;\;\text{on}\;\;\partial M\;\;\;\text{and}\;\;u(0)=0\;\;\text{on}\;\;M.
$$
Now given a function\;$v$\;on \;$\partial M$, we extend \;$v$\;(and denote still the extension by \;$v$)\;to a solution 
\begin{equation}\label{eq:tosolv}
 \left\{
\begin{split}
P^4_{g_0}w&=0\;\;&\text{on}\;\;M;\\
w&=v\;\;\;&\text{on}\;\;\partial M;\\
\frac{\partial w}{\partial n_{g_0}}&=0\;\;\;&\text{on}\;\;\partial M.
\end{split}
\right.
\end{equation}
Next, we define an operator on functions defined on \;$\partial M$\;as follows
$$
Av=-e^{-3u}P^3_{g_0}v.
$$
With this, we have that \;$u(t)$\;solves\;$\eqref{eq:tosolv}$,is equivalent to 
$$
\frac{\partial u(t)}{\partial t}=Au-(e^{-3u}T_{g_{0}}-\frac{\ov {T_{g(t)}}}{\ov S}S).
$$
We will study the latter problem in order to prove short-time existence for the flow under study. To do so, we will use linearization techniques. Hence it is useful to study the corresponding linearized equation. Therefore we start with the following Lemma.
\begin{lem}
 Let \;$f\in L^2(\partial M\times[0,1])$, and let us assume that \;$||u||_{L^{\infty}(\partial M\times[0,1])}\leq C$. Then if \;$v$\;solves
\begin{equation}
 \left\{
\begin{split}
 \frac{\partial v}{\partial t}&=Av+f\\
v(0)&=0
\end{split}
\right.
\end{equation}
then there holds
$$
||v||_{W^{1,2}(\partial M\times[0,1])}\leq C||f||_{L^{2}(\partial M\times[0,1])}.
$$
\end{lem}
\begin{pf}
Using the identity\;$(a-b)^2+2ab=a^2+b^2$, we get
$$
\int_0^1\int_{\partial M}(\frac{\partial v}{\partial t})^2dS_{g(t)}dt+\int_0^1\int_{\partial M}(Av)^2dS_{g(t)}dt=\int_0^1(\frac{\partial v}{\partial t}-Av)^2dS_{g(t)}dt+2\int_0^1\int_{\partial M}\frac{\partial v}{\partial t}AvdS_{g(t)}dt.
$$
On the other hand, using the expression of the operator\;$A$\;and the relation of the volume, we get
$$
\int_0^1\int_{\partial M}\frac{\partial v}{\partial t}AvdS_{g(t)}dt=-\int_0^1\int_{\partial M}\frac{\partial v}{\partial t}P^3_{g_0}vdS_{g_0}dt.
$$
Recalling that \;$P^4_{g_0}v=0$, we arrive to
$$
\int_0^1\int_{\partial M}\frac{\partial v}{\partial t}AvdS_{g(t)}dt=-\int_0^1<\frac{\partial v}{\partial t},P^{4,3}_{g_0}v>_{L^2(M,g_0)}dt.
$$
Now from the self-adjointness of \;$P^{4,3}_{g_0}$, we derive
$$
\int_0^1\int_{\partial M}\frac{\partial v}{\partial t}AvdS_{g(t)}dt=-\frac{1}{2}\int_0^1\frac{d}{dt}< v,P^{4,3}_{g_0}v>_{L^2(M,g_0)}dt.
$$
Thus using the fact that \;$v(0)=0$\;and \;$P^{4,3}_{g_0}$\;is non-negative, we get
$$
\int_0^1\int_{\partial M}\frac{\partial v}{\partial t}AvdS_{g(t)}dt\leq 0.
$$
Hence, we obtain
$$
\int_0^1\int_{\partial M}(\frac{\partial v}{\partial t})^2dS_{g(t)}dt+\int_0^1\int_{\partial M}(Av)^2dS_{g(t)}dt\leq\int_0^1(\frac{\partial v}{\partial t}-Av)^2dS_{g(t)}dt.
$$
Next, using the evolution equation solved by \;$v$, we get
$$
\int_0^1\int_{\partial M}(\frac{\partial v}{\partial t})^2dS_{g(t)}dt+\int_0^1\int_{\partial M}(Av)^2dS_{g(t)}dt\leq\int_0^1f^2dS_{g(t)}dt.
$$
Recalling that \;$u(t)$\;is bounded in time and space and \;$dS_{g(t)}=e^{3u(t)}dS_{g_0}$, we obtain
$$
\int_0^1\int_{\partial M}(\frac{\partial v}{\partial t})^2dS_{g_0}dt+\int_0^1\int_{\partial M}(P^3_{g_0}v)^2dS_{g_0}dt\leq\int_0^1f^2dS_{g_0}dt.
$$
Thus using the regularity result in Lemma\;$\ref{eq:reg}$, we derive
$$
||u||_{W^{1,2}(\partial M\times[0,1])}\leq C||f||_{L^2(\partial M\times[0,1])};
$$
as desired. Hence the Lemma is proved.
\end{pf}

The next Lemma is a higher-order analogue of the previous one. Precisely, we have
\begin{lem}\label{eq:injectivity}
 Let \;$m\in \N$, and\;$f\in W^{m,2}(\partial M\times[0,1])$. Assuming that \;$u$\;is bounded in \;$L^{\infty}(\partial M\times[0,1])$, in \;$W^{m-1,4}(\partial M\times[0,1])$, and in \;$W^{m-1,2}(\partial M\times[0,1])$, we have, if \;$v$\;is solution of the initial linear evolution equation
\begin{equation}
 \left\{
\begin{split}
 \frac{\partial v}{\partial t}&=Av+f\\
v(0)&=0
\end{split}
\right.
\end{equation}
then
$$
||v||_{W^{m-1,2}(\partial M\times[0,1])}\leq C||f||_{W^{m,2}(\partial  M\times[0,1])}.
$$
\end{lem}
To prove the Lemma, we will argue by induction on \;$m$. The Lemma above ensure the validity for \;$m=0$. Now let us suppose the it holds for \;$m$. By the assumptions on, we have
$$
||e^{3u}Av||_{W^{m-2,2}(\partial M\times[0,1])}\leq C||Av||_{W^{m-2,2}(\partial M\times[0,1])}.
$$
On the other hand, using the expression of \;$A$, and the equation solved by \;$v$, it is easily seen that \;$Av$\;solves the following initial value problem
\begin{equation}
 \left\{
\begin{split}
 \frac{\partial Av}{\partial t}&=A(Av)+Af-3\frac{\partial u}{\partial t}Av\\
v(0)&=0
\end{split}
\right.
\end{equation}
Therefore the hypothesis of induction imply
$$
||Av||_{W^{m-1,2}(\partial M\times[0,1])}\leq C||\frac{\partial u(t)}{\partial t}Av||_{W^{m-2,2}(\partial M\times[0,1])}+||Af||_{W^{m-2,2}(\partial M\times[0,1])}.
$$
Now, using the assumptions on \;$u(t)$, we get
$$
||Av||_{W^{m-1,2}(\partial M\times[0,1])}\leq C||P^3_{g_0}v||_{W^{m-2,4}(\partial M\times[0,1])}+||Af||_{W^{m-2,2}(\partial M\times[0,1])}.
$$
Next, by interpolation we obatin
$$
||Av||_{W^{m-1,2}(\partial M\times[0,1])}\leq C||P^3_{g_0}v||^{\frac{1}{2}}_{W^{m-2,2}(\partial M\times[0,1])}||P^3_{g_0}v||^{\frac{1}{2}}_{W^{m-1,2}(\partial M\times[0,1])}+||Af||_{W^{m-2,2}(\partial M\times[0,1])}.
$$
Hence we arrive to
$$
||Av||_{W^{m-1,2}(\partial M\times[0,1])}\leq C||v||^{\frac{1}{2}}_{W^{m+1,2}(\partial M\times[0,1])}||P^3_{g_0}v||^{\frac{1}{2}}_{W^{m-1,2}(\partial M\times[0,1])}+||f||_{W^{m+1,2}(\partial M\times[0,1])}.
$$
Using again the hyphothesis of induction, we get
$$
||Av||_{W^{m-1,2}(\partial M\times[0,1])}\leq C||f||^{\frac{1}{2}}_{W^{m,2}(\partial M\times[0,1])}||P^3_{g_0}v||^{\frac{1}{2}}_{W^{m-1,2}(\partial M\times[0,1])}+||f||_{W^{m+1,2}(\partial M\times[0,1])}.
$$
Thus we derive
$$
||P^3_{g_0}v||_{W^{m-1,2}(\partial M\times[0,1])}\leq C||f||^{\frac{1}{2}}_{W^{m,2}(\partial M\times[0,1])}||P^3_{g_0}v||^{\frac{1}{2}}_{W^{m-1,2}(\partial M\times[0,1])}+||f||_{W^{m+1,2}(\partial M\times[0,1])}.
$$
Hence we obtain
$$
||v||_{W^{m+2,2}(\partial M\times[0,1])}\leq C||f||_{W^{m+1,2}(\partial M\times[0,1])},
$$
as desired. This ends the proof of the Lemma.
\begin{pro}
Under the assumptions fo Theorem\;$\ref{eq:ttf}$, we have that the initial boundary value problem corresponding to \;$\eqref{eq:tflowbs}$, has a unique short-time solution \;$g(t)=e^{2u(t)}g_0$.
\end{pro}
\begin{pf}
 Writting \;$g(t)=e^{2u(t)}g_0$, we have that \;$g(t)$\;solves the initial boundary value problem corresponding to \;$\eqref{eq:tflowbs}$\;is equivalent to
\begin{equation}
 \left\{
\begin{split}
 \frac{\partial u}{\partial t}&=Au-(e^{-3u}T_{g_{0}}-\frac{\ov {T_{g(t)}}}{\ov S}S)\\
u(0)&=0
\end{split}
\right.
\end{equation}
Therefore our strategy will be to solve thie scalar evolution equation. To do so, we will look for zero of a operator. We define \;$G:W^{m+1.2}(\partial M\times[0,1])\rightarrow W^{m,2}(\partial M\times[0,1])$\;as follows
$$
G(u)=\frac{\partial u}{\partial t}-Au+(e^{-3u}T_{g_{0}}-\frac{\ov {T_{g(t)}}}{\ov S}S).
$$
Now we choose \;$u_0$\;such that \;$\frac{\partial^iG(u_0)}{\partial t^i}=0\;\;\forall\;1\leq i\leq m$, and \;$u_0$\;is bounded in \;$L^{\infty}(\partial M\times[0,1])$, in \;$W^{m-1,4}(\partial M\times[0,1])$, and in \;$W^{m,2}(\partial M\times[0,1])$. We have that the Frechet derivative of \;$G$\;at \;$u_0$\;is
$$
DG(u_0)w=\frac{\partial w}{\partial t}-Aw-3e^{-3u_0}w.
$$
Thus Lemma\;$\ref{eq:injectivity}$\; implies that the Linearization of \;$G$\;at \;$u_0$\;is bijective. Hence the Local Inversion theorem ensures that \;$G$\;is bijective around \;$u_0$.
On the other hand, since \;$u_0$\; satisfies \;$\frac{\partial^iG(u_0)}{\partial t^i}=0\;\;\forall\;1\leq i\leq m$, then there exists a function \;$w$\;closed to \;$G(u_0)$\;in the strong topology of \;$W^{m,2}(\partial M\times[0,1])$\; such that
\begin{equation}\label{eq:zero}
w(t)=0,\;\;\;t\in[0,\epsilon],
\end{equation}
for some positive and small \;$\epsilon$. Thus using the local invertibility of \;$G$\;around \;$u_0$, we get
$$
u=G^{-1}(w)
$$
is well-defined. Thus, from \;$\eqref{eq:zero}$, we infer that \;$u$\;is a short-time solution to our initial evolution problem, thus we have the existence. The uniqueness is consequence of Local inversion Theorem. Hence the Proposition is proved
\end{pf}
\footnotetext[1]{E-mail addresses:  ndiaye@sissa.it}


\begin{thebibliography}{99}


\bibitem{bran2} Branson T.P., {\em Differential operators canonically associated to a conformal structure}, Math. scand., 57-2 (1995), 293-345.

\bibitem{bo} Branson T.P., Oersted., {\em Explicit functional determinants in four dimensions}, Proc. Amer. Math. Soc 113-3(1991), 669-682.

\bibitem{bcy} Branson T.P., Chang S.Y.A., Yang P.C., {\em Estimates and extremal problems for the log-determinant on 4-manifolds}, Comm. Math. Phys., 149(1992), 241-262.

\bibitem{bren} Brendle S.,{\em Global existence and convergence for a higher order flow in conformal geometry}, Ann. of Math. 158 (2003),323-343.

\bibitem{br1} Brendle S., {\em Curvature flows on surfaces with boundary}, Math. Ann. 324, 491-519 (2002).

\bibitem{br2} Brendle S., {\em A family of curvature flows on surfaces with boundary}, Math. Z. 241, 829-869 (2002).

\bibitem{cgy1} Chang S.Y.A., Gursky M.J., Yang P.C., {\em An equation of Monge-Ampere type in conformal geometry, and four-manifolds of positive Ricci curvature}, Ann. of Math. 155-3(2002), 709-787.

\bibitem{cgy2} Chang S.Y.A., Gursky M.J., Yang P.C.,{\em A conformally invariant sphere theorem in four dimensions}, Publ. Math. Inst. Hautes etudes Sci. 9892003), 105-143.

\bibitem{cq1} Chang S.Y.A., Qing J.,{\em The Zeta Functional Determinants on manifolds with boundary 1. The Formula}, Journal of Functional Analysis 147, 327-362 (1997)

\bibitem{cq2} Chang S.Y.A., Qing J.,{\em The Zeta Functional Determinants on manifolds with boundary II. Extremal Metrics and Compactness of Isospectral Set}, Journal of Functional Analysis 147, 363-399 (1997)

\bibitem{cqy1} Chang S.Y.A., Qing J.,., Yang P.C.,{\em Compactification of a class of conformally flat 4-manifold}, Invent. Math. 142-1(2000), 65-93.

\bibitem{cqy2} Chang S.Y.A., Qing J.,., Yang P.C.,{\em On the Chern-Gauss-Bonnet integral for conformal metrics on\;${\bf R}^{4}$},\;Duke Math. J. 103-3(2000),523-544.

\bibitem{cy}  Chang S.Y.A.,  Yang P.C.,{\em Extremal metrics of zeta functional determinants on 4-manifolds}, ann. of Math. 142(1995), 171-212.

\bibitem{cy1}  Chang S.Y.A.,  Yang P.C.,{\em On a fourth order curvature invariant}.

\bibitem{dm1} Djadli Z., Malchiodi A., {\em A fourth order uniformization
theorem on some four manifolds with large total $Q$-curvature},
C.R.A.S., 340 (2005), 341-346.

\bibitem{ch} Chow B., {\em the Ricci flow on the 2-sphere}, J. Differential Geometry. 33 (1991), n0 2, 325-334.

\bibitem{fg} Fefferman C., Graham C,R., {\em Q-curvature and Poincar\'e metrics}, Mathematical Research Letters 9, 139-151(2002).

\bibitem{fg1} Fefferman C., Graham C., {\em Conformal invariants}, In Elie Cartan et les mathematiques d'aujourd'hui. Asterisque (1985), 95-116.

\bibitem{fh} Fefferman C., Hirachi Kengo., {\em Ambient metric construction of Q-curvature in conformal and CR geomotries}, Mathematical Research Letters .10, 819-831(2003).

\bibitem{fr} Friedman A., {\em Partial differential Equations}, Robert E. Krieger Publishing Company Malabar, Florida (1983).

\bibitem{gjms} Graham C,R., Jenne R., Mason L., Sparling G., {\em Conformally invariant powers of the laplacian, I:existence},  J.London Math.Soc 46(1992), no.2, 557-565.

\bibitem{gz} Graham C,R.,Zworsky M., {\em Scattering matrix in conformal geometry}. Invent math. 152,89-118(2003).

\bibitem{ha} Hamilton R., {\em The Ricci flow on surfaces}, Mathematics and general relativity (Santa Cruz, CA, 1986), 237-262, Contemp. Math., 71, Amer. Math. Soc., Providence, RI, !988.


\bibitem{ms} Malchiodi A., Struwe M., {\em $Q$-curvature flow on \;$S^4$\;}, J. Diff. Geom., 73-1 (2006),
 1-44.

\bibitem{nd} Ndiaye C.B., {\em Constant Q-curvature metrics in arbitrary dimension}
preprint 2006.

\bibitem{nd1} Ndiaye C.B., {\em Conformal metrics with constant \;$Q$-curvature for manifolds with boundary}, preprint SISSA.


\bibitem{nd2} Ndiaye C.B., {\em Constant \;$T$-curvature conformal metric on 4-manifolds with boundary}, preprint SISSA


\bibitem{p1} Paneitz S., {\em A quartic conformally covariant differential operator for arbitrary pseudo-Riemannian manifolds}, preprint, 1983.

\bibitem{p2} Paneitz S., {\em Essential unitarization of symplectics and applications to field quantization}, J. Funct. Anal. 48-3 (1982), 310-359.





%
%





\bibitem{si} Simon L., {\em Asymptotics for a class of non-linear evolution equations, with applications to geometry}, Ann of Math, 118 (1983), 525-571.
\end{thebibliography}
\end{document}